\documentclass[12pt]{article}

\usepackage{amsmath,amssymb,amsthm}
\usepackage{hyperref}
\usepackage{color}
\usepackage{bbm}

\newcommand{\ds}{\displaystyle}
\newcommand{\fr}[2]{\frac{#1}{#2}}
\newcommand{\dfr}[2]{\dfrac{#1}{#2}}
\newcommand{\cd}{\cdot}
\newcommand{\cds}{\cdots}
\newcommand{\dsum}{\displaystyle \sum}

\newcommand{\ul}[1]{\underline{#1}}

\renewcommand{\l}{\left}
\renewcommand{\r}{\right}

\newcommand{\vsb}{\vspace{2mm}}
\newcommand{\q}{\quad}

\newcommand{\la}{\langle}
\newcommand{\ra}{\rangle}

\newcommand{\abs}[1]{\lvert{#1}\rvert}

\newcommand{\maru}[1]{{\ooalign{\hfil#1\/\hfil\crcr
\raise.167ex\hbox{\mathhexbox20D}}}}

\newcommand{\Z}{\mathbb{Z}}
\newcommand{\C}{\mathbb{C}}

\newcommand{\N}{\mathbb{N}}
\newcommand{\Q}{\mathbb{Q}}

\DeclareMathOperator{\res}{Res}

\DeclareMathOperator{\aut}{Aut}
\DeclareMathOperator{\wt}{wt}
\DeclareMathOperator{\tr}{tr}

\DeclareMathOperator{\Hom}{Hom}
\DeclareMathOperator{\Sym}{Sym}

\DeclareMathOperator*{\tensor}{\otimes}

\newcommand{\vir}{\mathrm{Vir}}

\newcommand{\sfr}[2]{\leavevmode\kern-.1em
  \raise.5ex\hbox{\the\scriptfont0 #1}\kern-.1em
  /\kern-.15em\lower.25ex\hbox{\the\scriptfont0 #2}}

\newcommand{\shf}{\sfr{1}{2}}
\newcommand{\pii}{\pi \sqrt{-1}\, }
\newcommand{\hf}{\fr{1}{2}}

\newcommand{\w}{\omega}
\newcommand{\vacuum}{\mathbbm{1}}
\newcommand{\vac}{\vacuum}

\makeatletter
\@addtoreset{equation}{section}

\theoremstyle{plain}
\newtheorem{thm}{Theorem}[section]
\newtheorem{prop}[thm]{Proposition}
\newtheorem{lem}[thm]{Lemma}
\newtheorem{cor}[thm]{Corollary}

\theoremstyle{definition}
\newtheorem{df}[thm]{Definition}

\newtheorem{cond}{Condition}

\theoremstyle{remark}
\newtheorem{rem}[thm]{Remark}

\newcommand{\pf}{\noindent {\bf Proof:}\q }

\title{Extended Griess algebras and Matsuo-Norton trace formulae}

\author{Hiroshi Yamauchi\\
{\sl \small Department of Mathematics,  Tokyo Woman's Christian University}
\\
{\sl \small 2-6-1 Zempukuji, Suginami, Tokyo 167-8585 Japan}}

\date{}

\newcommand{\bmid}{~\Big|~}
\renewcommand{\ker}{\mathrm{Ker}\,}
\renewcommand{\o}{\mathrm{o}}
\newcommand{\longto}{\longrightarrow}
\newcommand{\VB}{{V\!\!B}}
\newcommand{\B}{\mathbb{B}}
\newcommand{\h}{\mathfrak{h}}

\begin{document}

\maketitle

\begin{abstract}
We introduce the $\Z_2$-extended Griess algebra of a vertex operator superalgebra
with an involution and derive the Matsuo-Norton trace formulae for the extended
Griess algebra based on conformal design structure.
We illustrate an application of our formulae by reformulating the one-to-one
correspondence between 2A-elements of the Baby-monster simple group and 
$N=1$ $c=7/10$ Virasoro subalgebras inside the Baby-monster vertex operator superalgebra.
\end{abstract}

\baselineskip 6mm

\tableofcontents

\section{Introduction}
A mysterious connection is known to exist between vertex operator algebras (VOAs) 
and finite simple groups.
One can explain that the $j$-invariant is made of the characters of the Monster 
simple group as a consequence of the modular invariance of characters of vertex 
operator algebras \cite{FLM,Zh}.
Matsuo \cite{Ma} introduced the notion {\it of class $\mathcal{S}^n$} of a VOA and
derived the formulae, which we will call the {\it Matsuo-Norton trace formulae}, 
describing trace of adjoint actions of the Griess algebra of a vertex operator algebra.
A VOA $V$ is called of class $\mathcal{S}^n$ if the invariant subalgebra 
under its automorphism group coincides with the subalgebra generated by 
the conformal vector up to degree $n$ subspace.
In the derivation of the formulae the non-associativity of 
of products of vertex operator algebras are efficiently used, resulting that 
the Matsuo-Norton trace formulae strongly encode structures of higher subspaces 
of vertex operator algebras.
Suitably applying the formula it outputs some information of structures 
related to automorphisms.
Here we exhibit an application of the formulae.

Let $V$ be a VOA of OZ-type with central charge $c$ and 
$e$ a simple $c=1/2$ Virasoro vector of $V$.
Then $V$ is a direct sum of irreducible $\la e\ra$-modules and the zero-mode
$\o(e)=e_{(1)}$ acts semisimply on $V$. 
It is shown in \cite{Mi1} that $\tau_e=\exp\l( 16\pii \o(e)\r)$ defines an element 
in $\aut(V)$.
The possible $\o(e)$-eigenvalues on the Griess algebra are 2, 0, 1/2 and 1/16.
(The eigenspace with eigenvalue 2 is 1-dimensional spanned by $e$.)
We denote by $d_\lambda$ the dimension of $\o(e)$-eigensubspace of $V_2$ with 
eigenvalue $\lambda$.
Then
\begin{equation}\label{eq:1.1}
  \tr_{V_2}\o(e)^i=2^i+0^i\, d_0+\l(\dfr{1}{2}\r)^i d_{1/2}
  +\l(\dfr{1}{16}\r)^i d_{1/16}, 
\end{equation}
where we have set $0^0=1$.
If $V$ is of class $\mathcal{S}^4$ then Matsuo-Norton trace formulae give
the values $\tr_{V_2}\o(e)^i$ for $i=0,1,2$ and if $c(5c+22)\ne 0$ we can 
solve \eqref{eq:1.1} and obtain 
\begin{equation}\label{eq:1.2}
\begin{array}{l}
  d_0=\dfr{(5c^2-100c+1188)\dim V_2-545c^2-2006c}{c(5c+22)},
  \vsb\\
  d_{1/2}=\dfr{-2((3c-110)\dim V_2+50c^2+192c)}{c(5c+22)},
  \vsb\\
  d_{1/16}=\dfr{64((2c-22)\dim V_2+10c^2+37c)}{c(5c+22)}.
\end{array}
\end{equation}
Therefore, we get an explicit form of the trace 
$\tr_{V_2}\tau_e=1+d_0+d_{1/2}-d_{1/16}$ as
\begin{equation}\label{eq:1.3}
  \tr_{V_2}\tau_e
  =\dfr{(5c^2-234c+2816)\dim V_2-1280c^2-4736c}{c(5c+22)}.
\end{equation}
If we apply this formula to the moonshine VOA $V^\natural$ \cite{FLM} 
then we see $\tr_{V^\natural_2}\o(e)=4372$ and $\tau_e$ belongs to 
the 2A-conjugacy class of the Monster (cf.~\cite{ATLAS}).
This example supports an existence of a link between the structure theory of VOAs and 
the character theory of finite groups acting on VOAs.

Partially motivated by Matsuo's work, H\"{o}hn introduced the notion of 
{\it conformal designs} based on vertex operator algebras in \cite{Ho4}, 
which would be a counterpart of block designs and spherical 
designs in the theories of codes and lattices, respectively.
The defining condition of conformal designs is already used in \cite{Ma} to derive
the trace formulae and H\"{o}hn sought and obtained results analogous to known ones
in block and spherical designs in \cite{Ho4}.
Contrary to the notion of of class $\mathcal{S}^n$, the definition of conformal 
design does not require the action of automorphism groups of VOAs.
Instead, it is formulated by the Virasoro algebra.
(The Virasoro algebra is the key symmetry in the 2-dimensional conformal field theory.)
The conformal design is a purely structure theoretical concept in the VOA theory 
and seems to measure a structural symmetry of VOAs, as one can deduce the Matsuo-Norton trace formulae 
via conformal designs.

The purpose of this paper is to extend Matsuo's work on trace formulae to 
vertex operator superalgebras (SVOAs) with an involution based on conformal design structure.
If the even part of the invariant subalgebra of an SVOA under an involution 
is of OZ-type then one can equip its weight 2 subspace with the structure of 
a commutative (but usually non-associative) algebra called the Griess algebra.
We extend this commutative algebra to a larger one by adding the odd part and call 
it the {\it $\Z_2$-extended Griess algebra}.
Our $\Z_2$-extended Griess algebra is still commutative but not super-commutative.
It is known that Virasoro vectors are nothing but idempotents in a Griess algebra.
In the odd part of the extended Griess algebra one can consider square roots of idempotents.
We will discuss the structure of the subalgebra generated by a square root of an idempotent
in the extended Griess algebra when the top weight of the odd part is small.

Then we derive trace formulae of adjoint actions on the odd part of the extended
Griess algebra based on conformal design structure.
Our formulae are a variation of Matsuo-Norton trace formulae.
As a main result of this paper we apply the trace formulae to the Baby-monster SVOA 
$\VB^\natural$ \cite{Ho1} and reformulate the one-to-one correspondence between the 2A-elements of 
the Baby-monster simple group and certain $c=7/10$ Virasoro vectors of $\VB^\natural$ 
obtained in \cite{HLY} in the supersymmetric setting by considering square roots 
of idempotents in the extended Griess algebra of $\VB^\natural$.
This result is in a sense suggestive.
It is shown in \cite{Ma} that if a VOA of OZ-type is of class $\mathcal{S}^8$ 
and its Griess algebra has dimension $d>1$ then the central charge is 24 and 
$d=196884$, those of the moonshine VOA.
For SVOAs it is shown in \cite{Ho4} that if the odd top level of an SVOA with top weight 3/2
forms a conformal 6-design then its central charge is either 16 or 47/2.
Using our formulae we will sharpen this result.
We will show if the odd top level of an SVOA with top weight 3/2
forms a conformal 6-design and in addition if it has a proper subalgebra isomorphic to
the $N=1$ $c=7/10$ Virasoro SVOA then its central charge is 47/2 and the odd top level is
of dimension 4371, those of the Baby-monster SVOA.
The author naively expects that the Baby-monster SVOA is the unique example 
subject to this condition.

The organization of this paper is as follows.
In Section 2 we review the notion of invariant bilinear forms on SVOAs.
Our definition of invariant bilinear forms on SVOAs is natural in the sense 
that if $M$ is a module over an SVOA $V$ then its restricted dual $M^*$ is also 
a $V$-module.
We also consider $\Z_2$-conjugation of invariant bilinear forms on SVOAs.
In Section 3 we introduce the $\Z_2$-extended Griess algebra for SVOAs with
involutions and consider square roots of idempotents.
When the odd top weight is less than 3 we describe possible structures of 
subalgebra generated by square roots under mild assumptions.
In Section 4 we derive trace formulae of adjoint actions on the odd part of 
the extended Griess algebra based on conformal design structure.
A relation between conformal design structure and generalized Casimir 
vectors is already discussed in \cite{Ma,Ho4} but we clarify it in our situation
to derive the formulae.
In Section 5 we apply our formulae to the VOAs with the odd top weight 1 
and to the Baby-monster SVOA.
In the latter case we reformulate the one-to-one correspondence in \cite{HLY}
in the supersymmetric setting.
Section 6 is the appendix and we list data of the trace formulae.

\paragraph{Acknowledgment.}
The author wishes to thank Professor Atsushi Matsuo for stimulating discussions 
and his Mathematica programs.
He also thanks Professor Masahiko Miyamoto for valuable comments.
Most of the results of this paper were obtained by using computer.
The author used a computer algebra system Risa/Asir for Windows.
This work was supported by JSPS Grant-in-Aid for Young Scientists (Start-up) 
No.~19840025 and (B) No.~21740011.

\paragraph{Notation and Terminology.}
In this paper we will work over the complex number field $\C$.
We use $\N$ to denote the set of non-negative integers.
A VOA $V$ is called {\it of OZ-type} if it has the $L(0)$-grading 
$V=\oplus_{n\geq 0}V_n$ such that $V_0=\C \vac$ and $V_1=0$.
(``OZ'' stands for ``Zero-One'', this is introduced by Griess.)
We denote the $\Z_2$-grading of an SVOA $V$ by $V=V^0\oplus V^1$.
We allow the case $V^1=0$.
For $a\in V^i$, we define its parity by $\abs{a}:=i \in \Z/2\Z$.
We assume that every SVOA in this paper has the $L(0)$-grading 
$V^0=\oplus_{n\geq 0} V_n$ and $V^1=\oplus_{n\geq 0} V_{n+k/2}$ with 
non-negative $k\in \Z$.
It is always assumed that $V_0=\C \vac$.
If $a\in V_n$ then we write $\wt(a)=n$.
If $V^1\ne 0$ then the minimum $h$ such that $V^1_h\ne 0$ is called 
{\it the top weight}  and the homogeneous subspace $V_h^1$ is called 
{\it the top level} of $V^1$.
A sub VOA $W$ of $V$ is called {\it full} if the conformal vector of $W$ is 
the same as that of $V$.
We denote the Verma module over the Virasoro algebra with central charge $c$ and 
highest weight $h$ by $M(c,h)$, and $L(c,h)$ denotes its simple quotient.
A Virasoro vector $e$ of $V$ with central charge $c_e$ is called {\it simple} 
if it generates a simple Virasoro sub VOA isomorphic to $L(c_e,0)$ in $V$.
Let $(M,Y_M(\cd,z))$ be a $V$-module and $g\in \aut(V)$.
The {\it $g$-conjugate} $M\!\circ g$ of $M$ is defined as $(M,Y^g_M(\cd,z))$ 
where $Y_M^g(a,z):=Y_M(ga,z)$ for $a\in V$. 
$M$ is called {\it $g$-stable} if $M\!\circ g \simeq M$, and {\it $G$-stable}
for $G<\aut(V)$ if $M$ is $g$-stable for all $g\in G$.
We write $Y(a,z)=\sum_{n\in \Z}a_{(n)}z^{-n-1}$ for $a\in V$ and define 
its {\it  zero-mode} by $\o(a):=a_{(\wt(a)-1)}$ if $a$ is homogeneous and extend linearly.
The supercommutator is denoted by $[\cd,\cd]_+$.

\section{Invariant bilinear form}

In this setion, we denote the complex number $e^{r\pii}$ by 
$\zeta^{r}$ for a rational number $r$.
In particular, $\zeta^n=(-1)^n$ for an integer $n$.
We denote the space of $V$-intertwining operators of type $M^1\times M^2 \to M^3$
by $\binom{M^3}{M^1~M^2}_V$.

\subsection{Self-dual module}

Let $V$ be a simple self-dual VOA and $M$ an irreducible self-dual $V$-module.
We assume that $M$ has the $L(0)$-decomposition $M=\oplus_{n=0}^\infty M_{n+h}$ 
such that each homogeneous component is finite dimensional and 
the top weight $h$ of $M$ is either in $\Z$ or in $\Z+1/2$.
Let $( \cd \mid \cd )_V$ and $( \cd \mid \cd )_M$ be invariant bilinear
forms on $V$ and $M$, respectively.
One can define $V$-intertwining operators $I(\cd,z)$ and $J(\cd,z)$ of types 
$M\times V \to M$ and $M\times M \to V$, respectively, as follows 
(cf.~Theorem 5.5.1 of \cite{FHL}).
For $a \in V$ and $u,v\in M$, 
\begin{equation}\label{eq:2.1}
\begin{array}{l}
  I(u,z)a:=e^{zL(-1)} Y_M(a,-z) u,
  \vsb\\
  ( J(u,z)v \mid a)_V
  := (u \mid I(e^{zL(1)} (\zeta z^{-2})^{L(0)} u,z^{-1})a)_M.
\end{array}
\end{equation}
(Note that $I(u,z)a\in M(\!(z)\!)$ and $J(u,z)v\in V(\!(z)\!)$.)
Since $\dim \binom{M}{M~V}_V=\dim \binom{V}{M~M}_V=1$, 
it follows from Proposition 2.8 of \cite{Li1} and Proposition 5.4.7 of \cite{FHL} 
that there exist $\alpha,\beta\in \{\pm 1\}$ such that
\begin{equation}\label{eq:2.2}
  ( u \mid v )_M = \alpha ( v \mid u)_M,\q
  e^{zL(-1)} J(v,-z) u = \beta J(u,z)v \q 
  \text{for}\  u,v\in M.
\end{equation}
We can sharpen Proposition 5.6.1 of \cite{FHL} as follows.

\begin{lem}\label{lem:2.1}
  Let $\alpha,\beta\in \{\pm 1\}$ de defined as above.
  Then $\alpha\beta =(-1)^{2h}$.
\end{lem}

\pf
Let $a\in V$ and $u,v\in M$.
$$
\begin{array}{ll}
  ( J(u,z)v \mid a )_V 
  &= \l( v \mid I(e^{zL(1)} (\zeta z^{-2})^{L(0)} u,z^{-1})a\r)_M
  \vsb\\
  &= \l( v \bmid e^{z^{-1}L(-1)} Y_M(a,-z^{-1}) e^{zL(1)} (\zeta z^{-2})^{L(0)} u \r)_M
  \vsb\\
  &= \l( e^{z^{-1}L(1)} v \bmid Y_M(a,-z^{-1}) e^{zL(1)} (\zeta z^{-2})^{L(0)} u \r)_M
  \vsb\\
  &= \l( Y_M(e^{-z^{-1} L(1)} (-z^2)^{L(0)} a,-z) e^{z^{-1}L(1)} v \bmid 
     e^{zL(1)} (\zeta z^{-2})^{L(0)} u \r)_M
  \vsb\\
  &= \l( e^{-zL(-1)} I(e^{z^{-1}L(1)} v,z) e^{-z^{-1} L(1)} (-z^2)^{L(0)} a \bmid
     e^{zL(1)} (\zeta z^{-2})^{L(0)} u \r)_M
  \vsb\\
  &= \l( I(e^{z^{-1}L(1)} v,z) e^{-z^{-1} L(1)} (-z^2)^{L(0)} a \bmid
      (\zeta z^{-2})^{L(0)} u \r)_M
  \vsb\\
  &= \alpha \l( (\zeta z^{-2})^{L(0)} u \bmid I(e^{z^{-1}L(1)} v,z) 
     e^{-z^{-1} L(1)} (-z^2)^{L(0)} a \r)_M .
\end{array}
$$
By the definition of the invariance, one has 
$$
  \l( v \mid I(u,z) a\r)_M 
  = \l( J(e^{zL(1)}(\zeta^{-1}z^{-2})^{L(0)} u,z^{-1}) v \mid a \r)_V.
$$
Then we continue 
$$
\begin{array}{l} 
  \alpha \l( (\zeta z^{-2})^{L(0)} u \bmid I(e^{z^{-1}L(1)} v,z) 
     e^{-z^{-1} L(1)} (-z^2)^{L(0)} a \r)_M 
  \vsb\\
  = \alpha \Big( J(e^{zL(1)} 
  \underbrace{(\zeta^{-1} z^{-2})^{L(0)} e^{z^{-1}L(1)}}_{=e^{-zL(1)}(\zeta^{-1} z^{-2})^{L(0)}}
  v,z^{-1}) (\zeta z^{-2})^{L(0)} u \bmid  e^{-z^{-1} L(1)} (-z^2)^{L(0)} a \Big)_V
  \vsb\\
  = \alpha \l( J((\zeta^{-1} z^{-2})^{L(0)} v,z^{-1})  
     (\zeta z^{-2})^{L(0)} u \bmid  e^{-z^{-1} L(1)} (-z^2)^{L(0)} a \r)_V
  \vsb\\
  = \alpha \l( e^{-z^{-1} L(-1)}J((\zeta^{-1} z^{-2})^{L(0)} v,z^{-1})  
     (\zeta z^{-2})^{L(0)} u \bmid (-z^2)^{L(0)} a \r)_V
  \vsb\\
  = \alpha \beta \l( J((\zeta z^{-2})^{L(0)} u,-z^{-1})(\zeta^{-1} z^{-2})^{L(0)} v
        \bmid (-z^2)^{L(0)} a \r)_V
  \vsb\\
  = \alpha \beta \underbrace{\zeta^{\wt(u)-\wt(v)}}_{=(-1)^{\wt(u)-\wt(v)}} 
    \l( J(z^{-2L(0)} u,-z^{-1})z^{-2L(0)} v \mid (-z^2)^{L(0)} a \r)_V
  \vsb\\
  = \alpha \beta (-1)^{\wt(u)-\wt(v)}\l((-z^2)^{L(0)}J(z^{-2L(0)} u,-z^{-1})
    z^{-2L(0)} v \mid  a \r)_V.
  \vsb\\
\end{array}
$$
Since 
$(-z^2)^{L(0)}J(u,-z^{-1})v=(-1)^{\wt(u)+\wt(v)}J(z^{2L(0)}u,z)z^{2L(0)}v$, 
we further continue 
$$
\begin{array}{l}
  \alpha \beta \zeta^{\wt(u)-\wt(v)}((-z^2)^{L(0)}J(z^{-2L(0)} u,-z^{-1})z^{-2L(0)} v
        \mid  a )_V
  \vsb\\
  = \alpha \beta (-1)^{\wt(u)-\wt(v)}\cd (-1)^{\wt(u)+\wt(v)} (J(u,z)v \mid  a )_V
  \vsb\\
  = \alpha \beta (-1)^{2\wt(u)} (J(u,z)v \mid  a )_V
  \vsb\\
  = \alpha \beta (-1)^{2h} ( J(u,z)v \mid a )_V.
\end{array}
$$
Therefore, we obtain the desired equation $\alpha \beta =(-1)^{2h}$.
\qed 

\subsection{Invariant bilinear forms on SVOA}

Let $V=V^0\oplus V^1$ be an SVOA.
Let $M$ be an untwisted $V$-module, that is, $M$ has a $\Z_2$-grading 
$M=M^0\oplus M^1$ compatible with that of $V$.
We also assume that $M^i$ ($i=0,1$) has an $L(0)$-decomposition 
$M^i=\oplus_{n\geq 0}M^i_{h+n+i/2}$ where $h$ is the top weight of $M$ and 
each $L(0)$-eigensubspace is finite dimensional.
Let $M^*$ be its restricted dual.
We can define a vertex operator map $Y_{M^*} (\cd ,z)$ on $M^*$ by means of 
the adjoint action $\l\la Y_{M^*} (a,z)\nu \mid v \r\ra = \la \nu \mid Y_M^*(a,z)v \ra$
for $a \in V$, $\nu \in M^*$ and $v \in M$, where
\begin{equation}\label{eq:2.3}
  Y^*_M(a,z):=Y_M(e^{zL(1)} z^{-2L(0)}\zeta^{L(0)+2L(0)^2} a,z^{-1}).
\end{equation}
Then one can show the structure $(M^*,Y^*_M (\cd,z))$ forms a $V$-module (cf.~\cite{FHL}).

\begin{lem}\label{lem:2.2}
  $(M^*,Y_{M^*}(\cd,z))$ is a $V$-module.
\end{lem}

\begin{rem}
  The correction term $\zeta^{2L(0)^2}$ in the definition of the adjoint operator
  $Y^*_M(\cd,z)$ is necessary for $M^*$ to be a $V$-module.
  It follows from our assumption on the $L(0)$-grading that 
  $\zeta^{L(0)+2L(0)^2}a =\pm a$ if $a\in V$ is $\Z_2$-homogeneous.
  Therefore, $\zeta^{-L(0)-2L(0)^2}a =\zeta^{L(0)+2L(0)^2}a$ and we have 
  $Y^{**}_M(a,z)=Y_M(a,z)$.
\end{rem}

If $M^*$ is isomorphic to $M$ as a $V$-module, then $M$ is called {\it self-dual} 
(as a $V$-module) and there exists an invariant bilinear form 
$( \cd \mid \cd )_M$ on $M$ satisfying
\begin{equation}\label{eq:2.4}
  \l( Y(a,z)u \mid v \r)_M 
  = \l( u \bmid Y(e^{zL(1)} z^{-2L(0)}\zeta^{L(0)+2L(0)^2} a,z^{-1}) v\r)_M
\end{equation}
for $a\in V$ and $u,v\in M$.

The following is an easy generalization of known results (cf.~Lemma \ref{lem:2.1}).

\begin{prop}(\cite{FHL, Li1})
  Let $V$ be an SVOA of CFT-type.
  \\
  (1) Any invariant bilinear form on $V$ is symmetric.
  \\
  (2) The space space of invariant bilinear forms on $V$ is linearly isomorphic 
  to the dual space of $V_0/L(1)V_1$, i.e., $\Hom_\C (V_0/L(1)V_1,\C)$.
  In particular, if $V^1$ is irreducible over $V^0$ then 
  $V$ is a self-dual $V$-module if and only if $V^0$ is a self-dual $V^0$-module.
\end{prop}

\subsection{Conjugation of bilinear forms}


Let $\theta=(-1)^{2L(0)}$ be the canonical $\Z_2$-symmetry of a superalgebra, i.e., 
$\theta=1$ on $V^0$ and $\theta=-1$ on $V^1$, then we have 
\begin{equation}\label{eq:2.5}
  Y^*_M(\theta a,z)=(-1)^{2\wt(a)}Y^*_M(a,z)
  =Y_M(e^{zL(1)} z^{-2 L(0)}\zeta^{L(0)-2L(0)^2} a,z^{-1}).
\end{equation}
Consider the $\theta$-conjugate $M\!\circ\theta =(M, Y^\theta_M(\cd,z))$ of $M$.
If $M=M^0\oplus M^1$ is also a superspace then $M$ is always $\theta$-stable, for, 
one can define an isomorphism $\tilde{\theta}$ between $M$ and $M\circ \theta$ by 
$\tilde{\theta}=1$ on $M^0$ and $\tilde{\theta}=-1$ on $M^1$.
Therefore, if $M$ is self-dual then $M$, $M^*$ and $(M\!\circ \theta)^*\simeq M^*\!\circ \theta$ are all isomorphic.
This means there is a freedom of choice of the adjoint operator in the case of SVOA.
We can choose the right hand side of \eqref{eq:2.5} as well as \eqref{eq:2.3} for 
the definition of $Y^*(\cd,z)$.
From now on we will freely choose one of \eqref{eq:2.3} or 
\eqref{eq:2.5} for the adjoint operator.

\section{Extended Griess algebras}

In this section we introduce a notion of $\Z_2$-extended Griess algebras.

\subsection{Definition}

We will consider SVOAs subject to the following condition.

\begin{cond}\label{cond:1}
Let $V=V^0\oplus V^1$ be a vertex operator superalgebra of CFT-type 
and $g$ an involution of $V$.
Denote $V^\pm :=\{ a\in V \mid \pm g a=a\}$.
We assume the following.
\\
(1)~$V$ is self-dual.
\\
(2)~$V$ has the $L(0)$-grading $V^0=\oplus_{n\geq 0}V_n$ and 
$V^1=\oplus_{n\geq 0}V_{n+k/2}$ with non-negative $k\in \Z$ (if $V^1$ is non-zero).
\\
(3)~$V^\pm$ has the $L(0)$-decomposition such that
$V^+=V_0\oplus V_2 \oplus \l(\oplus_{n>2} V_n\r)$ and 
$V^- = V_h \oplus \l( \oplus_{n>h} V_h\r)$ where $V_h\ne 0$ is the
top level of $V^-$ and $h\in \hf\Z$ is its top weight.
\end{cond}

Note that $V_h\subset V^0$ if $h\in \Z$ and $V_h\subset V^1$ if $h\in \Z+1/2$.
We will denote $V^{0,+}:=V^0\cap V^+$.
By assumption, $V^{0,+}$ is of OZ-type and $V_2$ is the Griess algebra of $V^{0,+}$.

If $V^1\ne 0$ and $h\in \Z+1/2$ we choose
\begin{equation}\label{eq:3.1}
  Y^*(a,z)=
  \begin{cases}
    Y(e^{zL(1)} z^{-2L(0)}\zeta^{L(0)-2L(0)^2} a,z^{-1}) & \mbox{ if } h \equiv 1/2 \mod 2,
    \vsb\\
    Y(e^{zL(1)} z^{-2L(0)}\zeta^{L(0)+2L(0)^2} a,z^{-1}) & \mbox{ if } h \equiv 3/2 \mod 2,
\end{cases}
\end{equation}
so that if we write $Y^*(a,z)=\dsum_{n\in \Z}a^*_{(n)}z^{-n-1}$ then 
\begin{equation}\label{eq:3.2}
  a^*_{(n)}= \epsilon_h (-1)^{\wt(a)-h} \dsum_{i=0}^\infty \dfr{1}{i!}
  \l( L(1)^i a\r)_{(2\wt(a)-n-2-i)} ,
\end{equation}
where the signature $\epsilon_h\in \{ \pm 1\}$ is defined by
\begin{equation}\label{eq:3.3}
  \epsilon_h=\begin{cases} (-1)^h & \mbox{ if } h\in \Z,
  \vsb\\ ~~1 & \mbox{ if } h\in \Z+1/2. \end{cases}
\end{equation}
In particular, we have $(u|v)\vac = \epsilon_h u_{(2h-1)}v$ for $u,v\in V_h$.

\medskip

Consider the subspace $V_2\oplus V_h$ of $V^+\oplus V^-$.

\begin{prop}
  For $a,b\in V_2$ and $u,v\in V_h$, define 
  \begin{equation}\label{eq:3.4}
  \begin{array}{l}
  a b:=a_{(1)}b,~~~
  a u:=a_{(1)}u,~~~
  u a:=u_{(1)}a,~~~
  u v:=u_{(2h-3)}v,
  \vsb\\
  (a|b)\vac=a_{(3)}b,~~~
  (a|u)=0=(u|a),~~~
  (u|v)\vac=\epsilon_h  u_{(2h-1)}v.
  \end{array}
  \end{equation}
  Then the subspace $V_2\oplus V_h$ forms a unital 
  commutative $\Z_2$-graded algebra with invariant bilinear form 
  which extends the Griess algebra structure on $V_2$, where the invariance 
  of the bilinear form is modified as $(xu|y)= \epsilon_h (x|uy)$ for $u\in V_h$.
\end{prop}

\pf
The proof follows by a direct verification.
For example, by the skew-symmetry one has
$$
\begin{array}{ll}
  uv
  &\ds =u_{(2h-3)}v
  =\dsum_{j\geq 0} \dfr{(-1)^{2h-3+1+j+\abs{u}\cd \abs{v}}}{j!}L(-1)^j v_{(2h-3+j)}u
  \vsb\\
  &\ds = \underbrace{(-1)^{2h+\abs{u}\cd \abs{v}}}_{=\, 1} \l( v_{(2h-3)}u
  +(-1)\cd L(-1)\underbrace{v_{(2h-2)}u}_{\in\, V_1^+\, =\, 0}
  +\dfr{(-1)^2}{2}\underbrace{L(-1)^2v_{(2h-1)}u}_{\in \, L(-1)^2V_0=0}\r)
  \vsb\\
  &= vu.
\end{array}
$$
That $au=ua$ also follows from the skew-symmetry.
The bilinear form clearly satisfies the invariant property.
\qed
\vsb

We call $V_2\oplus V_h$ the {\it $\Z_2$-extended Griess algebra} of $V$.

\subsection{Square roots of idempotents}\label{sec:3.1}

Recall the following fact (cf.~Lemma 5.1 of \cite{Mi1} and Proposition 2.6 of \cite{La}).

\begin{prop}
  Let $V$ be a VOA of CFT-type.
  A vector $e\in V_2$ is a Virasoro vector with central charge $c$ 
  if and only if it satisfies $e_{(1)}e=2e$ and $2e_{(3)}e=c\vac$.
\end{prop}

By this proposition, for a VOA $V$ of OZ-type we see that $e\in V_2$ is a 
Virasoro vector if and only if $e/2$ is an idempotent in the Griess algebra.
So idempotents are important objects to study in the Griess algebra.
If we consider the $\Z_2$-extended Griess algebra, it is possible to consider 
square roots of idempotents inside the odd part. 

Let $V=V^0\oplus V^1$ be an SVOA and $\theta=(-1)^{2L(0)} \in \aut(V)$ the canonical 
$\Z_2$-symmetry of $V$.
In this subsection we assume $V$ and $g=\theta$ satisfy Condition \ref{cond:1} and 
we consider its extended Griess algebra $V_2\oplus V_h$.
Let $a$ be an idempotent of $V_2$ and suppose $x\in V_h$ is a square root of $a$, 
that is, $xx=a$ hold in the extended Griess algebra.
We shall consider the subalgebra $\la x\ra$ generated by such a root $x$.
The structure of $\la x\ra$ depends on the top weight $h$.

\paragraph{Case $h=1/2$.}
In this case $x_{(n)}x\in V_{-n}$ and $x_{(n)}x=0$ if $n>0$ as $V$ is of CFT-type.
Since $x_{(0)}x=(x|x)\vac$, we have the following commutation relation:
\begin{equation}\label{eq:3.5}
  [x_{(m)},x_{(n)}]_+
  =\l( x_{(0)}x\r)_{(m+n)}=(x|x)\vac_{(m+n)}=(x|x)\delta_{m+n,-1}.
\end{equation}
Since $a=xx=x_{(-2)}x$ is an idempotent, we have $aa=a_{(1)}a=a$.
Then 
$$
\begin{array}{ll}
  (x_{(-2)}x)_{(1)}(x_{(-2)}x)
  &= \dsum_{i\geq 0}(-1)^i\binom{-2}{i} \l(x_{(-2-i)}x_{(1+i)}
     -(-1)^{-2+\abs{x}\cd\abs{x}} x_{(-1-i)}x_{(i)}\r) x_{(-2)}x
  \vsb\\
  &= x_{(-2)}x_{(1)}x_{(-2)}x+x_{(-1)}x_{(0)}x_{(-2)}x+2x_{(-2)}x_{(1)}x_{(-2)}x
  \vsb\\
  &= 4(x|x) x_{(-2)}x
\end{array}
$$
and therefore we have $4(x|x)=1$.
The central charge of $2a$ is $8(a|a)=8(x_{(-2)}x|x_{(-2)}x)=8(x|x_{(1)}x_{(-2)}x)
=8(x|(x|x)\cd x)=1/2$.
Set $\psi_{n+1/2}:=2x_{(n)}$ and $\psi(z):=Y(x,z)$. 
Then \eqref{eq:3.5} is expressed as $[\psi_r,\psi_s]_+=\delta_{r+s,0}$ for 
$r,s\in \Z+1/2$ and we have a free fermionic field 
$$
  \psi(z)=\dsum_{r\in \Z+1/2} \psi_r z^{-r-1/2},~~~~
  \psi(z)\psi(w)\sim \dfr{1}{z-w}.
$$
Therefore, $\la x\ra$ is isomorphic to a simple $c=1/2$ Virasoro SVOA 
$L(\shf,0)\oplus L(\shf,\shf)$ (cf.~\cite{KR}).

\paragraph{Case $h=3/2$.}
In this case we have $x_{(0)}x=a$, $x_{(1)}x=0$, $x_{(2)}x=(x|x)\vac$ and 
$x_{(n)}x=0$ for $n\geq 3$.
Here we further assume that the $c=8(a|a)$ Virasoro vector $2a$ is 
the conformal vector of the subalgebra $\la x\ra$. 
This condition is equivalent to $x\in \ker_V (\w -2a)_{(0)}$ (cf.~\cite{FZ}), 
where $\w$ is the conformal vector of $V$.
(By Theorem 5.1 of (loc.~cit), $\w-2a$ is a Virasoro vector and $\w=2a+(\w-2a)$ 
is an orthogonal decomposition.)
Then $x$ is a highest weight vector for $\vir(2a)$ with highest weight $3/2$. 
Since $xx=x_{(0)}x=a$, we have $(a|x_{(0)}x)=(a|a)$.
On the other hand, 
$$
\begin{array}{ll}
  (2a|x_{(0)}x)\vac &= 2a_{(3)}x_{(0)}x=[2a_{(3)},x_{(0)}]x
  = \dsum_{i\geq 0} \binom{3}{i}\l( 2a_{(i)}x\r)_{(3-i)}x
  \vsb\\
  &= \l( 2a_{(0)}x\r)_{(3)}x+3\l( 2a_{(1)}x\r)_{(2)} 
    = -3 x_{(2)}x+3\cd \dfr{3}{2}x_{(2)}x
  \vsb\\
  &= \dfr{3}{2}(x|x)\vac
\end{array}
$$
and we get $3(x|x)=4(a|a)$.
Now we can compute the commutation relations:
\begin{equation}\label{eq:3.6}
\begin{array}{ll}
\begin{array}{l}
  [2a_{(m)},x_{(n)}]
  = \dsum_{i=0}^\infty \binom{m}{i} (2a_{(i)}x)_{(m+n-i)}
  \vsb\\
  = (2a_{(0)}x)_{(m+n)}+m(2a_{(1)}x)_{(m+n-1)}
  \vsb\\
  = -(m+n)x_{(m+n-1)}+m\cd \dfr{3}{2}x_{(m+n-1)}
  \vsb\\
  = \dfr{1}{2}(m-2n)x_{(m+n-1)}, 
\end{array}
&
\begin{array}{l}
  {}[x_{(p)},x_{(q)}]_+
  = \dsum_{i=0}^\infty \binom{p}{i} (x_{(i)}x)_{(p+q-i)}
  \vsb\\
  = \ds (x_{(0)}x)_{(p+q)}+ \binom{p}{2} (x_{(2)}x)_{(p+q-2)}
  \vsb\\
  = a_{(p+q)}+\dfr{p(p-1)}{2} (x|x) \vac_{(p+q-2)}
  \vsb\\
  = a_{(p+q)}+\dfr{2p(p-1)}{3} (a|a) \delta_{p+q,1}.
\end{array}
\end{array}
\end{equation}
For simplicity, set $L^a(m):=2a_{(m+1)}$, $c_a:=8(a|a)$ and $G^x(r):=2x_{(r+1/2)}$.
Then \eqref{eq:3.6} looks 
\begin{equation}\label{eq:3.7}
\begin{array}{l}
  [L^a(m), G^x(r)]=\dfr{1}{2}(m-2r) G^x(m+r),
  \vsb\\
  {}[G^x(r),G^x(s)]_+=2L^a(r+s)+\delta_{r+s,0}\dfr{4r^2-1}{12} c_a .
\end{array}
\end{equation}
This is exactly the defining relations of the Neveu-Schwarz algebra, also known as 
the $N=1$ super Virasoro algebra.
Therefore, $\la x\ra$ is isomorphic to the $N=1$ $c=8(a|a)$ Virasoro SVOA.

\paragraph{Case $h=5/2$.}
Again we assume that $xx=a$ and $2a$ is the conformal vector of $\la x\ra$.
But this is still not enough to determine the structure of $\la x\ra$.
In order to describe $\la x\ra$, we need one more assumption that 
$x_{(n)}x\in \la a\ra$ for $n\geq 0$.
Zamolodchikov \cite{Za} has studied such a subalgebra.

\begin{prop}\label{prop:3.3}(\cite{Za})
  Suppose $x_{(n)}x\in \la a\ra$ for $n\geq 0$ and $2a$ is the conformal vector 
  of $\la x\ra$.
  Then there is a surjection from $\la x\ra$ to $L(-\sfr{13}{14},0)\oplus 
  L(-\sfr{13}{14},\sfr{5}{2})$.
  In particular, the central charge of $\la x\ra$ is uniquely determined.
\end{prop}

\begin{rem}
Recall the central charges and the highest weights of the the minimal series 
of the Virasoro VOAs (cf.~\cite{W}). 
\begin{equation}\label{eq:3.8}
  c_{p,q}:=1-\dfr{6(p-q)^2}{pq},~~
  h^{(p,q)}_{r,s}=\dfr{(rq-sp)^2-(p-q)^2}{4pq},~~
  0<r<p,~ 0<s<q.
\end{equation} 
Then $c_{7,4}=-13/14$ and $h^{(7,4)}_{6,1}=5/2$.
Moreover, $L(-\sfr{13}{14},\sfr{5}{2})$ is the unique non-trivial 
simple current module over $L(-\sfr{13}{14},0)$ so that the simple quotient
in Proposition \ref{prop:3.3} forms a $\Z_2$-graded simple current extension
of $L(-\sfr{13}{14},0)$.
\end{rem}

\begin{rem}
The (extended) Griess algebra is a part of the structure of the vertex Lie algebra 
\cite{P} (or that equivalently known as the conformal algebra \cite{K}).
As seen in this subsection, for small $h$ one can determine OPE of elements 
in $V_2\oplus V_h$ by the extended Griess algebra.
However, for higher $h$, the extended Griess algebra is insufficient to determine 
full OPE of $Y(x,z)$ and the subalgebra $\la x\ra$ for $x\in V_2\oplus V_h$. 
\end{rem}

\section{Matsuo-Norton Trace formulae for extended Griess algebras}

In this section we derive a trace formulae for the extended Griess algebras, 
which is a variation of Matsuo-Norton trace formulae \cite{Ma}. 
In this section we assume $V$ satisfies the the following condition.

\begin{cond}\label{cond:2}
  $V$ is an SVOA satisfying Condition \ref{cond:1} and in addition the following.
  \\
  (1)~The invariant bilinear form is non-degenerate on $V$.
  \\
  (2)~The restriction of the bilinear form on $\vir(\w)$ is also non-degenerate.
  \\
  (3)~$V$ as a $\vir(\w)$-module is a direct sum of highest weight modules. 
\end{cond}

Let $V[n]$ be the sum of highest weight $\vir(\w)$-submodules of $V$ with highest 
weight $n \in \hf \Z$.
Then by (3) of Condition \ref{cond:2} one has 
\begin{equation}\label{eq:4.1}
  V=\oplus_{n\geq 0} V[n]
\end{equation}
and we can define the projection map
\begin{equation}\label{eq:4.2}
  \pi : V=\bigoplus_{n\geq 0} V[n]\longto V[0]=\vir(\w)=\la \w\ra
\end{equation}
which is a $\vir(\w)$-homomorphism.

\subsection{Conformal design and Casimir vector}

Let us recall the notion of the conformal design.
Suppose $U$ is a VOA satisfying (3) of Condition \ref{cond:2}.
Then we can define the projection $\pi : U \to \la \w\ra$ as in \eqref{eq:4.2}.

\begin{df}(\cite{Ho4})
  Let $U$ be a VOA and $M$ a $U$-module.
  An $L(0)$-homogeneous subspace $X$ of $M$ is called a {\it conformal $t$-design} based on $V$ 
  if $\tr_X\o(a) =\tr_X\o(\pi(a))$ holds for any $a\in \oplus_{0\leq n\leq t} U_n$. 
\end{df}

The defining condition of conformal designs was initiated in Matsuo's paper \cite{Ma} 
and it is related to the following condition.

\begin{df}(\cite{Ma})
  Let $U$ be a VOA and $G$ a subgroup of $\aut(U)$.
  We say $U$ is {\it of class $\mathcal{S}^n$} under $G$ if 
  $U^G_k \subset \la \w\ra$ for $0\leq k\leq n$. 
  (We allow $G$ to be $\aut(U)$ itself.)
\end{df}

The above two conditions are in the following relation.

\begin{lem}(\cite{Ma,Ho4})
  Let $U$ be a VOA and $M$ a $G$-stable $U$-module.
  Suppose further that we have a projective representation of $G$ on $M$.
  If $U$ is of class $\mathcal{S}^n$ under $G$ then every $L(0)$-homogeneous 
  subspace of $M$ forms a conformal $n$-design.
\end{lem}

Let $V$ be an SVOA and $g\in \aut(V)$ satisfying Condition \ref{cond:2}.
The following is clear. 

\begin{lem}\label{lem:4.4}
  For $m>0$ the components $V[m]$ and $V[0]$ in \eqref{eq:4.1} are orthogonal 
  with respect to the invariant bilinear form.
\end{lem}

Consider the extended Griess algebra $V_2\oplus V_h$ of $V$.
Set $d:=\dim V_h$.
Let $\{ u_i\}_{1\leq i\leq d}$ be a basis of $V_h$ and $\{ u^i\}_{1\leq i\leq d}$ 
its dual basis.
As in \cite{Ma}, we consider the {\it Casimir vector} of weight $m\in \Z$: 
\begin{equation}\label{eq:casimir}
  \kappa_m:= \epsilon_h \dsum_{i=1}^d u^i_{(2h-1-m)} u_i \in V_m,
\end{equation}
where the signature $\epsilon_h$ is defined as in \eqref{eq:3.3}.

\begin{lem}\label{lem:4.5}
  Let $a\in V$ be homogeneous.
  Then $\tr_{V_h}\o(a)=(-1)^{\wt(a)} (a\mid \kappa_{\wt(a)})$.
\end{lem}

\pf
Without loss we may assume $a$ is even, i.e., $a\in V^{0,+}$.
We compute the trace as follows.
$$
\begin{array}{lll}
  \tr_{V_h}\o(a) 
  &= \dsum_{i=1}^d (\o(a) u^i | u_i) =\dsum_{i=1}^d (a_{(\wt(a)-1)} u^i \mid u_i) & 
  \vsb\\
  &= \dsum_{i=1}^d \dsum_{j=0}^\infty \dfr{(-1)^{\wt(a)+j}}{j!}
  \l( L(-1)^j u^i_{(\wt(a)-1+j)}a \bmid u_i\r)
  & \mbox{(by skew-symmetry)}
  \vsb\\
  &= \dsum_{i=1}^d \dsum_{j=0}^\infty \dfr{(-1)^{\wt(a)+j}}{j!} \epsilon_h 
  \l( a \bmid u^i_{(2h-\wt(a)-1+j)}L(1)^j u_i\r)& \mbox{(by invariance)}
  \vsb\\
  &= \dsum_{i=1}^d (-1)^{\wt(a)}\epsilon_h (a|u^i_{(2h-\wt(a)-1)}u_i) 
  & \mbox{(as }L(1)V_h=0)
  \vsb\\
  &= (-1)^{\wt(a)}(a \mid \kappa_{\wt(a)}) .&
\end{array}
$$
Therefore, we obtain the desired equality.
\qed

\begin{prop}\label{prop:4.6}
  $V_h$ is a conformal $t$-design based on $V^{0,+}$ 
  if and only if $\kappa_m\in \la \w\ra$ for $0\leq m\leq t$.
\end{prop}

\pf
By Lemma \ref{lem:4.5}, $\tr_{V_h}\o(a)=\tr_{V_h}\o(\pi(a))$ for any $a\in V^{0,+}_m$
if and only if $(a| \kappa_m) =(\pi(a) | \kappa_m)$.
Since $\pi$ is a projection, $\{ a-\pi(a) \mid a \in V^{0,+}_m\} = V^{0,+}_m \cap \ker \pi$.
Then $(a-\pi(a) | \kappa_m)=0$ if and only if $\kappa_m \in \pi(V)=V[0]=\la \w\ra$ 
by Condition \ref{cond:2} and Lemma \ref{lem:4.4}. 
Therefore the assertion holds.
\qed

\subsection{Derivation of trace formulae}
We use Lemma \ref{lem:4.5} to derive the trace formulae.
Recall the following associativity formula.

\begin{lem}\label{lem:4.7}
(\cite[Lemma 3.12]{Li2})
Let $a$, $b\in V$, $v\in M$ and $p$, $q\in \Z$.
Suppose $s\in \Z$ and $t \in \N$ satisfy
$a_{(s+i)}v=b_{(q+t+i+1)}v=0$ for all $i\geq 0$.
Then for $p$, $q\in \Z$, 
\begin{equation}\label{eq:4.4}
  a_{(p)}b_{(q)}v=\dsum_{i=0}^t \dsum_{j\geq 0}\binom{p-s}{i}\binom{s}{j}
  \l( a_{(p-s-i+j)}b\r)_{(q+s+i-j)}v.
\end{equation}
\end{lem}

Let $a, b\in V_2$ and $v\in V_h$.
Since $V_h$ is the top level of $V^-$ and $V_2\subset V^+$, 
we can apply the lemma above with $p=q=1$, $s=2$ and $t=0$ and obtain 
\begin{equation}\label{eq:4.5}
  \o(a)\o(b)v=a_{(1)}b_{(1)}v=  \dsum_{i=0}^2\binom{2}{i}\l( a_{(i-1)}b\r)_{(3-i)}
  =\o(a*b),
\end{equation}
where $a*b$ is the product of Zhu algebra \cite{Zh} defined as 
\begin{equation}\label{eq:4.6}
  a*b:=\res_z \dfr{(1+z)^{\wt (a)}}{z}Y(a,z)b
  =\dsum_{i\geq 0} \binom{\wt(a)}{i} a_{(i-1)}b.
\end{equation}

Combining this with Lemma \ref{lem:4.5}, we obtain

\begin{lem}\label{lem:4.8}
  For $a^1,\dots,a^k \in V_2$ we have 
  \begin{equation}\label{eq:4.7}
  \begin{array}{l}
    \tr_{V_h}\o(a^1)\cds \o(a^k)
    \vsb\\
    =\dsum_{-1\leq i_1,\cds,i_{k-1}\leq 1}
    (-1)^{k+1-(i_1+\cds +i_{k-1})} 
    \l(a^1_{(i_1)}\cds a^{k-1}_{(i_{k-1})} a^k \bmid 
    \kappa_{k+1-(i_1+\cds +i_{k-1})}\r) .
  \end{array}
  \end{equation}
\end{lem}

\pf
By \eqref{eq:4.5} if $a\in V_2$ and $x\in V_n$ then 
$a*x=a_{(1)}x+2a_{(0)}x+a_{(-1)}x \in V_n+V_{n+1}+V_{n+2}$.
So by Lemma \ref{lem:4.5}
\begin{equation}\label{eq:4.8}
  \tr_{V_h}\o(a^1)\cds \o(a^k)=\dsum_{m=2}^{2k} (-1)^m 
  \l(\o(a^1\!\!*\! \cds\! *\!a^k)\mid \kappa_m\r) .
\end{equation}
Expanding this we obtain the lemma.
\qed
\vsb

To describe the Casimir vector we need the following condition.

\begin{cond}\label{cond:3}
  If $V_h$ forms a conformal $2t$-design with $t\leq 5$
  then the central charge of $V$ is not a zero of the polynomial 
  $D_{2t}(c)$ defined as follows.
  \begin{equation}\label{eq:4.9}
  \begin{array}{l}
     D_2(c)=c,~~~
     D_4(c)=c(5c+22),~~~
     D_6(c)=(2c-1)(7c+68)D_4(c),
     \vsb\\
     D_8(c)=(3c+46)(5c+3)D_6(c),~~~~
     D_{10}(c)=(11c+232)D_8(c).
   \end{array}
  \end{equation}
\end{cond}

The normalized polynomials $D_n(c)$ comes from the Shapovalov determinant of 
the Verma module $M(c,0)$.
The following is well-known (cf.~\cite{KR}).

\begin{lem}
  If the central charge of $V$ is not a zero of $D_n(c)$ in \eqref{eq:4.9} 
  then the degree $m$ subspace of $\la \w\ra$ with $m\leq n$ is isomorphic to 
  that of $M(c,0)/M(c,1)$.
\end{lem}

We write $[n_1,\dots,n_k]\models m$ if $n_1\geq \cds \geq n_k\geq 2$ and $n_1+\cds +n_k=m$.
If $c$ is not a zero of $D_n(c)$ in \eqref{eq:4.9} then the degree $m$ subspace 
of $M(c,0)/M(c,1)$ with $m\leq n$ has a basis 
$\{ L(-n_1)\cds L(-n_k)\vac \mid [n_1,\dots,n_k] \models m\}$.
By definition $\kappa_0=\epsilon_h \sum_{i=1}^d u^i_{(2h-1)}u_i 
= \sum_{i=1}^d (u^i|u_i)\vac = d\vac$ and $\kappa_1=0$, where $d=\dim V_h$.
For $n>0$ one can show
\begin{equation}\label{eq:4.10}
  L(n)\kappa_m=\l( h(n-1)+m-n\r) \kappa_{m-n}.
\end{equation}
Using this we can recursively compute $(L(-n_1)\cds L(-n_k)\mid \kappa_m)$.
As a result, the Casimir vector can be expressed as follows (cf.~Proposition 2.5 of \cite{Ma}).

\begin{lem}\label{lem:4.10}
  Suppose the central charge of $V$ is not a zero of $D_n(c)$ in \eqref{eq:4.9} and 
  the Casimir vector $\kappa_m\in \la \w\ra$ with $m\leq n$.
  Then $\kappa_m$ is uniquely written as
  $$
    \kappa_m = \dfr{1}{D_{2\lfloor m/2\rfloor}(c)}\dsum_{[n_1,\dots,n_k] \,\models\, m} 
    A^{(m)}_{[n_1,\dots,n_k]} L(-n_1)\cds L(-n_k)\vac, 
  $$
  where $\lfloor m/2\rfloor$ stands for the largest integer not exceeding $m/2$ and
  $A^{(m)}_{[n_1,\dots,n_k]} \in \Q[c,d,h]$ are given in Appendix \ref{app:casimir}. 
\end{lem}

We present the main result of this paper.

\begin{thm}\label{thm:4.11}
  Suppose $V$ and $g\in \aut(V)$ satisfy Conditions \ref{cond:2} and \ref{cond:3}.
  Set $d=\dim V_h$.
  \\
  (1)~If $V_h$ forms a conformal 2-design based on $V^{0,+}$, then for any $a^0\in V_2$, 
  $$
    \tr_{V_h}\o(a^0)= \dfr{2hd}{c}(a^0|\w).
  $$
  (2)~If $V_h$ forms a conformal 4-design based on $V^{0,+}$, then for any $a^0,a^1\in V_2$, 
  $$
    \tr_{V_h}\o(a^0)\o(a^1)= \dfr{4hd(5h+1)}{c(5c+22)} (a^0|\w)(a^1|\w)
    + \dfr{2hd(22h-c)}{c(5c+22)} (a^0|a^1).
  $$
  (3)~If $V_h$ forms a conformal 6-design based on $V^{0,+}$, then for any $a^0,a^1,a^2\in V_2$, 
  $$
  \begin{array}{l}
    \tr_{V_h}\o(a^0)\o(a^1)\o(a^2) 
    \vsb\\
    = D_6(c)^{-1} \l( F^{(3)}_0 (a^0|\w)(a^1|\w)(a^2|\w) 
    +F^{(3)}_1 \Sym (a^0|\w)(a^1|a^2) 
    +F^{(3)}_2 (a^0|a^1|a^2)\r),
  \end{array}
  $$
  where $(a^0|a^1|a^2)=(a^0a^1|a^2)$ is a totally symmetric trilinear form, 
  $\Sym (a^0|\w)(a^1|a^2)$ is the sum of all $(a^{i_0}|\w)(a^{i_1}|a^{i_2})$ which are 
  mutually distinct, and $F^{(3)}_j \in \Q[c,d,h]$, $0\leq j\leq 2$, are given in Appendix \ref{app:formula}.
  \\
  (4)~If $V_h$ forms a conformal 8-design based on $V^{0,+}$, then for any $a^0,a^1,a^2,a^3\in V_2$, 
  \\
  $$
  \begin{array}{l}
  \tr_{V_h}\o(a^0)\o(a^1)\o(a^2)\o(a^3) 
  \vsb\\
  = \ds  D_8(c)^{-1} \l( F^{(4)}_0 (a^0|\w)(a^1|\w)(a^2|\w)(a^3|\w)  
  +F^{(4)}_1 \Sym (a^0|\w)(a^1|\w)(a^2|a^3)
  \r.
    \vsb\\
    \ds ~~~+F^{(4)}_2 \Sym (a^0|\w)(a^1|a^2|a^3)
    +F^{(4)}_3 \Sym (a^0|a^1)(a^2|a^3)
        +F^{(4)}_4 (a^0a^1|a^2a^3)
    \vsb\\
    ~~~~ \l.+F^{(4)}_5 (a^0a^2|a^1a^3) +F^{(4)}_6 (a^0a^3|a^1a^2)\r) ,
  \end{array}
  $$
  where $\Sym$ denotes the sum over all possible permutations of $(0,1,2,3)$ 
  for which we obtain mutually distinct terms, 
  and $F^{(4)}_j \in \Q[c,d,h]$, 
  $0\leq j\leq 6$, are given in Appendix \ref{app:formula}.
  \\
  (5)~If $V_h$ forms a conformal 10-design based on $V^{0,+}$, then for any $a^0,a^1,a^2,a^3,a^4\in V_2$, 
  $$
  \begin{array}{l}
    \tr_{V_h}\o(a^0)\o(a^1)\o(a^2)\o(a^3)\o(a^4)
    \vsb\\
    = D_{10}(c)^{-1} \l( 
     F^{(5)}_0 (a^0|\w)(a^1|\w)(a^2|\w)(a^3|\w)(a^4|\w) 
     +F^{(5)}_1 \Sym (a^0|\w)(a^1|\w)(a^2|\w)(a^3|a^4)
     \r.
    \vsb\\ 
~~~~    +F^{(5)}_2\!\Sym (a^0|\w)(a^1|\w)(a^2|a^3|a^4)
    +F^{(5)}_3\!\Sym (a^0|\w)(a^1|a^2)(a^3|a^4)
    \vsb\\ 
~~~~    +F^{(5)}_4\big(
    (a^0|\w)(a^1a^2|a^3a^4)
    +(a^1|\w)(a^0a^2|a^3a^4)
    +(a^2|\w)(a^0a^1|a^3a^4)
    \vsb\\ 
~~~~~~~~~~~~~~~~~~     +(a^3|\w)(a^0a^1|a^2a^4)
     +(a^4|\w)(a^0a^1|a^2a^3)\big)
    \vsb\\ 
~~~~    +F^{(5)}_5\big(
    (a^0|\w)(a^1a^3|a^2a^4)
    +(a^1|\w)(a^0a^3|a^2a^4) 
    +(a^2|\w)(a^0a^3|a^1a^4)
    \vsb\\ 
~~~~~~~~~~~~~~~~~~    +(a^3|\w)(a^0a^2|a^1a^4)
    +(a^4|\w)(a^0a^2|a^1a^3)\big)
    \vsb\\ 
~~~~    +F^{(5)}_6\big(
    (a^0|\w)(a^1a^4|a^2a^3)
    +(a^1|\w)(a^0a^4|a^2a^3)
    +(a^2|\w)(a^0a^4|a^1a^3)
    \vsb\\ 
~~~~~~~~~~~~~~~~~~    +(a^3|\w)(a^0a^4|a^1a^2)
    +(a^4|\w)(a^0a^3|a^1a^2)\big)
    \vsb\\ 
~~~~+F^{(5)}_7\Sym (a^0|a^1)(a^2|a^3|a^4)
    +F^{(5)}_8 (a^0a^1a^2a^3a^4)
    \l. 
    +{\dsum}^{*} F^{(5)}_{i_0i_1i_2i_3i_4} (a^{i_0}a^{i_1}|a^{i_2}|a^{i_3}a^{i_4}) \r)
  \end{array}
  $$
  where $\Sym$ denotes the sum over all possible permutations of $(0,1,2,3,4)$ 
  for which we obtain mutually distinct terms, 
  $(a^0a^1a^2a^3a^4)\vac=a^0_{(3)}a^1_{(2)}a^2_{(1)}a^3_{(0)}a^4$ and 
  the last summation ${\sum}^*$ is taken over all possible permutations $(i_0,i_1,i_2,i_3,i_4)$ of 
  $(0,1,2,3,4)$ such that $(a^{i_0}a^{i_1}|a^{i_2}|a^{i_3}a^{i_4})$ are mutually distinct. 
  The coefficients $F^{(5)}_{\bullet}\in \Q[c,d,h]$ are given in Appendix \ref{app:formula}.
\end{thm}

\pf
By Lemmas \ref{lem:4.8} and \ref{lem:4.10}, it suffices to rewrite inner products 
of the form 
$(a^0_{(i_0)}\cds a^{k-1}_{(i_{k-1})}a^k \mid L(-m_1)\cds L(-m_l)\vac)$, 
$-1\leq i_0, \cds,i_{k-1}\leq 1$, $k+1-(i_0+\cds i_{k-1})=m_1+\cds +m_k$, 
in terms of the Griess algebra. 
By the invariance, this is equal to 
$(L(m_l)\cds L(m_1)a^0_{(i_0)}\cds a^{k-1}_{(i_{k-1})}a^k\mid \vac)$, 
and by the commutation formula
$$
  [L(m),a_{(n)}]=(m-n+1)a_{(m+n)}+\delta_{m+n,1}\dfr{m(m^2-1)}{6}(a|\w),
$$
we obtain a sum of $(a^{s_0}_{(j_0)}\cds a^{s_{r-1}}_{(j_{r-1})}a^{s_r} \mid \vac)$ with $r\leq k$ 
and $-1\leq j_0,\cds,j_{r-1}< 2k$.
We will use the following relations to rewrite such a term further.
\begin{equation}\label{eq:4.11}
\begin{array}{l}
  (a_{(0)}b)_{(n)}=[a_{(1)},b_{(n-1)}]-(a_{(1)}b)_{(n-1)},
  ~~~a_{(-m-1)}\vac = \dfr{1}{m!}L(-1)^m a,
  \vsb\\
  {}\big[ a_{(m)}\,,\,b_{(0)}\big] =\big[ a_{(1)}\,,\,b_{(m-1)}\big]+(m-1)(a_{(1)}b)_{(m-1)},
  ~~~a_{(m)}b_{(-m)}\vac = m a_{(1)}b, 
\end{array}
\end{equation}
where $a,b\in V_2$, $m>0$ and $n\in \Z$.
For example, let us rewrite 
$( a^0_{(4)}a^1_{(1)}a^2_{(0)}a^3 \mid \vac)$.
$$
\begin{array}{ll}
  \l( a^0_{(4)}a^1_{(1)}a^2_{(0)}a^3\bmid\vac\r)
  &= \l( a^2_{(0)}a^3 \bmid a^1_{(1)}a^0_{(-2)}\vac\r)
  =\l( a^2_{(0)}a^3 \bmid \big[ a^1_{(1)}\,,\, a^0_{(-2)}\big]\vac\r)
  \vsb\\
  &= \l( a^2_{(0)}a^3 \bmid \big( a^1_{(0)}a^0\big)_{(-1)}\vac\r)
  + \l( a^2_{(0)}a^3 \bmid \big( a^1_{(1)}a^0\big)_{(-2)}\vac\r)
  \vsb\\
  &= \l( a^2_{(0)}a^3 \bmid a^1_{(0)}a^0\r)
  + \l( a^3 \bmid a^2_{(2)} \big( a^1_{(1)}a^0\big)_{(-2)}\vac\r)
  \vsb\\
  &= \l( a^0_{(3)}a^1_{(2)}a^2_{(0)}a^3 \bmid \vac\r)
  + 2\l( a^3 \bmid a^2_{(1)}a^1_{(1)}a^0\r) 
  \vsb\\
  &= \l( a^0_{(3)}a^1_{(2)}a^2_{(0)}a^3 \bmid \vac\r)
  + 2\l( a^2_{(1)}a^3 \bmid a^1_{(1)}a^0\r) ,
\end{array}
$$
where $( a^0_{(3)}a^1_{(2)}a^2_{(0)}a^3|\vac)$ can be simplified as
$$
\begin{array}{ll}
  \l( a^0_{(3)}a^1_{(2)}a^2_{(0)}a^3 \bmid \vac\r)
  &= \l( a^1_{(2)}a^2_{(0)}a^3 \bmid a^0 \r)
  = \l( \big[ a^1_{(2)}\,,\, a^2_{(0)}\big] a^3 \bmid a^0 \r)
  \vsb\\
  &= \l( \big[ a^1_{(1)}\,,\, a^2_{(1)}\big] a^3 \bmid a^0 \r)
  + \l( \big( a^1_{(1)}a^2\big)_{(1)}a^3 \bmid a^0\r)
  \vsb\\
  &= \l( a^2_{(1)}a^3 \bmid a^1_{(1)}a^0\r) 
  -\l( a^1_{(1)}a^3 \bmid a^2_{(1)}a^0\r) 
  +\l( a^1_{(1)}a^2 \bmid a^3_{(1)}a^0\r) .
\end{array}
$$
Therefore, we get 
$$
  \l( a^0_{(4)}a^1_{(1)}a^2_{(0)}a^3\bmid\vac\r)
  = 3\l( a^0_{(1)}a^1 \bmid a^2_{(1)}a^3\r) 
  -\l( a^0_{(1)}a^2 \bmid a^1_{(1)}a^3\r) 
  +\l( a^0_{(1)}a^3 \bmid a^1_{(1)}a^2\r) .
$$
In this way, we can rewrite all terms and obtain the formulae of degrees 1 to 4.
However, in the rewriting procedure of the trace of degree 5 we meet the expressions
$a^{i_0}_{(3)}a^{i_1}_{(2)}a^{i_2}_{(1)}a^{i_3}_{(0)}a^{i_4}=(a^{i_0}a^{i_1}a^{i_2}a^{i_3}a^{i_4})\vac$ 
which satisfy the following relations:
\begin{equation}\label{eq:4.12}
\begin{array}{ll}
  (a^0a^1a^2a^3a^4) +(a^1a^0a^2a^3a^4)  
  &=
  3(a^0a^1|a^2|a^3a^4)
  -(a^0a^1|a^3|a^2a^4)
  +(a^0a^1|a^4|a^2a^3),
  \vsb\\
  (a^0a^1a^2a^3a^4)+(a^0a^2a^1a^3a^4)
  &=(a^0a^1|a^2|a^3a^4)-(a^0a^1|a^3|a^2a^4)+(a^0a^1|a^4|a^2a^3)
  \vsb\\
  &+(a^0a^2|a^1|a^3a^4)-(a^0a^2|a^3|a^1a^4)+(a^0a^2|a^4|a^1a^3)
  \vsb\\
  &-(a^0a^3|a^4|a^1a^2)+(a^0a^4|a^3|a^1a^2)+(a^1a^2|a^0|a^3a^4),
  \vsb\\
  (a^0a^1a^2a^3a^4)+(a^0a^1a^3a^2a^4)
  &=(a^0a^1|a^2|a^3a^4)+(a^0a^1|a^3|a^2a^4)+(a^0a^1|a^4|a^2a^3)
  \vsb\\
  &-(a^0a^2|a^1|a^3a^4)-(a^0a^3|a^1|a^2a^4)+(a^0a^4|a^1|a^2a^3)
  \vsb\\
  &+(a^1a^2|a^0|a^3a^4)+(a^1a^3|a^0|a^2a^4)-(a^1a^4|a^0|a^2a^3),
  \vsb\\
  (a^0a^1a^2a^3a^4)+(a^0a^1a^2a^4a^3)
  &=3(a^0a^1|a^2|a^3a^4)-(a^0a^2|a^1|a^3a^4)+(a^1a^2|a^0|a^3a^4).
\end{array}
\end{equation}
Let $R$ be the space of formal sums of $(a^{i_0}a^{i_1}|a^{i_2}|a^{i_3}a^{i_4})$
over $\Z$.
Then it follows from \eqref{eq:4.12} that
$$
  (a^{\sigma(0)}a^{\sigma(1)}a^{\sigma(2)}a^{\sigma(3)}a^{\sigma(4)})
  \equiv \mathrm{sign}(\sigma) (a^0a^1a^2a^3a^4) \mod R
$$
for $\sigma\in \mathrm{S}_5$.
This shows the rewriting procedure is not unique, and in our rewriting procedure 
we have to include at least one term $(a^0a^1a^2a^3a^4)$ in the formula (cf.~\cite{Ma}).
\qed

\begin{rem}\label{rem:4.12}
Our formulae are a variation of Matsuo-Norton trace formulae in \cite{Ma} but 
there are some differences.
The trace formula of degree 5 in (loc.~cit) contains a totally anti-symmetric quinary form 
$$
  \dsum_{\sigma \in \mathrm{S}_5} \mathrm{sgn}(\sigma) 
  \big( a^{\sigma(0)}a^{\sigma(1)}a^{\sigma(2)}a^{\sigma(3)}a^{\sigma(4)}\big)
  = \dsum_{\sigma \in \mathrm{S}_5} \mathrm{sgn}(\sigma) \,
  \l( a^{\sigma(0)}_{(3)}a^{\sigma(1)}_{(2)}a^{\sigma(2)}_{(1)}a^{\sigma(3)}_{(0)}a^{\sigma(4)}\bmid \vac\r) ,
$$
which we do not have in ours.
This is due to the non-uniqueness of the reduction procedure as explained in the proof.
One can transform the formula to include this form using \eqref{eq:4.12}.
\end{rem}

\begin{rem}\label{rem:4.13}
In the trace formula of degree $n$, if we put $a^i=\w/h$ for one of $0\leq i< n$ 
then we obtain the trace formula of degree $n-1$.
Even though we have derived the formulae for degree 4 and 5, the author does not know 
non-trivial examples of SVOAs which satisfy Conditions \ref{cond:2} and \ref{cond:3} and the 
odd top level forms a conformal 4- or 5-design.
It is shown in \cite{Ma} (see also \cite{Ho4} for related discussions) that 
if a VOA $V$ satisfying Conditions \ref{cond:2} and \ref{cond:3} and is 
of class $\mathcal{S}^8$ (under $\aut(V)$) and has a proper idempotent then 
$\dim V_2=196884$ and $c=24$, those of the moonshine VOA.
By this fact, the author expected the non-existence of proper SVOAs (not VOAs) 
of class $\mathcal{S}^8$ and $\mathcal{S}^{10}$, but the reductions of the trace formula 
of degree 5 to degree 4 and of degree 4 to degree 3 are consistent and 
we cannot obtain any contradiction.
\end{rem}

We will mainly use the formulae to compute traces of Virasoro vectors.

\begin{cor}\label{cor:4.14}
  Suppose $V$ and $g\in \aut(V)$ satisfy Conditions \ref{cond:2} and \ref{cond:3}.
  Let $e\in V_2$ be a Virasoro vector with central charge $c_e=2(e|e)$.
  If $V_h$ forms a conformal $2t$-design based on $V^{0,+}$, then $\tr_{V^h} \o(e)^t$ 
  is given as follows.
  \\
  $$
  \begin{array}{l}
    \tr_{V^h}\o(e) =\dfr{2hd}{c}(e|e)~~~ \mbox{if }~t=1,~
    \vsb\\
    \tr_{V^h}\o(e)^2=\dfr{4hd(5h+1)}{c(5c+22)}(e|e)^2+\dfr{2hd(22h-c)}{c(5c+22)}(e|e)~~~
    \mbox{if }~t=2,
    \vsb\\
    \tr_{V^h}\o(e)^t = D_{2t}(c)^{-1}\dsum_{j=1}^t E^{(t)}_j (e|e)^j ~~~\mbox{if }~t=3,4,5,
  \end{array}
  $$
  where $d=\dim V_h$ and $E^{(t)}_\bullet \in \Q[c,d,h]$ are given in Appendix \ref{app:formula}.
\end{cor}

\section{Applications}

In this section we show some applications of our formulae.

\subsection{VOAs with $h=1$}

Let $V$ be a VOA and $\theta\in \aut(V)$ an involution satisfying 
Conditions \ref{cond:1} and \ref{cond:2}.
We assume the top weight $h$ of $V^-$ is $1$ and denote $d=\dim V_1$.
In this case the top level $V_1$ forms an abelian Lie algebra under 0-th product.
For, $[V_1,V_1]=(V_1)_{(0)}V_1\subset V_1^+=0$ as $V_1\subset V^-$ and $V^+$ is of OZ-type.
We have the following commutation relation for $a$ and $b\in V_1$:
\begin{equation}\label{eq:5.1}
  [a_{(m)},b_{(n)}]=(a_{(0)}b)_{(m+n)}+m(a_{(1)}b)_{(m+n-1)}=-m(a|b)\delta_{m+n,0}.
\end{equation}
Set $\h :=V_1$ and equip this with a symmetric bilinear form by $\la a|b\ra:=-(a|b)$ 
for $a,b\in \h$.
Then this form is non-degenerate by Condition \ref{cond:2}.
Denote $\hat{\h}$ the rank $d=\dim V_1$ Heisenberg algebra associated to $\h$, 
the affinization of $\h$.
By \eqref{eq:5.1} the sub VOA $\la V_1\ra$ generated by $V_1$ is isomorphic to 
a free bosonic VOA associated to $\hat{\h}$ and the Casimir element $\kappa_2$ is 
twice the conformal vector of $\la V_1\ra$.
Suppose $V_1$ forms a conformal 2-design based on $V^+$.
Then $\kappa_2$ coincides with twice the conformal vector of $V$, and hence 
$\la V_1\ra$ is a full sub VOA of $V$.
More precisely, we show the following.

\begin{prop}\label{prop:5.1}
  Suppose a VOA $V$ and its involution $\theta\in \aut(V)$ satisfy Conditions 
  \ref{cond:1} and \ref{cond:2}.
  If the top level of $V^-$ has the top weight $1$ and forms a conformal 2-design 
  based on $V^+$ then $\la V_1\ra$ is isomorphic to the rank $d=\dim V_1$ free bosonic 
  VOA and the restriction of $\theta$ on $\la V_1\ra$ is conjugate to a lift of the 
  $(-1)$-map on $\h$ in $\aut(\la V_1\ra)$.
  Moreover, if $\la V_1\ra$ is a proper subalgebra of $V$ then there is an even
  positive definite rootless lattice $L$ of $\h$ of rank less than or equal 
  to $d$ such that $V$ is isomorphic to a tensor product of the lattice VOA $V_L$
  associated to $L$ and a free bosonic VOA associated to the affinization of 
  the orthogonal complement of $\C L$ in $\h$.
  In this case the restriction of $\theta$ on $V_L$ is conjugate to a lift of the 
  $(-1)$-map on $L$ in $\aut(V_L)$.
\end{prop}

\pf
That $\la V_1\ra$ is isomorphic to a free bosonic VOA is already shown and 
$\theta$ is clearly a lift of the $(-1)$-map on it.
Suppose $\la V_1\ra$ is a proper subalgebra.
Then by \cite{LX} there is a even positive definite lattice $L$ such that 
$V$ is isomorphic to a tensor product of the lattice VOA $V_L$ associated to $L$ 
and a free bosonic VOA associated to the affinization of the orthogonal 
complement of $\C L$ in $\h$.
Since $V_1$ is abelian, $L$ has no root.
Let $\rho$ be a lift of the $(-1)$-map on $L$ to $\aut(V_L)$.
Since $\theta$ and $\rho$ act by $-1$ on $(V_L)_1$, $\theta\rho$ is identity on it.
Then $\theta\rho$ on $V_L$ is a linear automorphism $\exp(a_{(0)})$ for 
some $a \in (V_L)_1$.
Since $\theta\exp(\hf a_{(0)})=\exp(\hf (\theta a)_{(0)})\theta=\exp(-\hf a_{(0)})\theta$, 
we have the conjugacy $\exp(-\hf a_{(0)}) \theta \exp (\hf a_{(0)})
=\theta \exp(\hf a_{(0)})\exp (\hf a_{(0)}) =\theta \exp (a_{(0)})
=\theta(\theta \rho)=\rho$.
This completes the proof.
\qed
\vsb

If $V$ is a free bosonic VOA then it is shown in \cite{DN} that 
$V^+$ is not of class $\mathcal{S}^4$.
The case $V=\la V_1\ra$ is not interesting and out of our focus. 
So we are reduced to the case when $V$ is a lattice VOA $V_L$ 
where $L$ is an even positive definite rootless lattice and 
$\theta$ is a lift of the $(-1)$-map on $L$.
(Clearly $V_L$ and $\theta$ satisfy Conditions \ref{cond:1} and \ref{cond:2}.)
For such a lattice $L$ the complete classification of simple $c=1/2$ Virasoro
vectors in $V_L^+$ is obtained in \cite{Sh2}.
It is shown in \cite{LSY,Sh2} that $V_L$ has a simple $c=1/2$ Virasoro vector 
$e\in V_L^+$ if and only if there is a sublattice $K$ isomorphic to
$\sqrt{2}A_1$ or $\sqrt{2}E_8$ such that $e\in V_K^+\subset V_L^+$.
In particular, one can find a simple $c=1/2$ Virasoro vector in $V_L$ 
if (and only if) $L$ has a norm 4 element.

Suppose we have a simple $c=1/2$ Virasoro vector $e\in V_L$.
Actually $e$ is in $V_L^+$ in our case and one can find a sublattice $K$
isomorphic to $\sqrt{2}A_1$ or $\sqrt{2}E_8$ such that $e\in V_K^+\subset V_L^+$.
Denote $X=(V_L)_1$.
The zero-mode $\o(e)$ acts on $X$ semisimply with possible eigenvalues 
$0$, $1/2$ and $1/16$.
We denote by $d_\lambda$ the dimension of $\o(e)$-eigensubspace of $X$ with eigenvalue $\lambda$.
Then $d_0+d_{1/2}+d_{1/16}=d$.
By the trace formula in (1) of Theorem \ref{thm:4.11} we get
\begin{equation}\label{eq:5.2}
  \tr_X \o(e)=0\cd d_0 + \dfr{1}{2}d_{1/2}+\dfr{1}{16}d_{1/16}=\dfr{1}{2}.
\end{equation}
This has two possible solutions $(d_0,d_{1/2},d_{1/16})=(d-1,1,0)$ and $(d-8,0,8)$.
Suppose further that $X$ forms a conformal 4-design.
Then by the trace formula in (2) of Theorem \ref{thm:4.11} we have
\begin{equation}\label{eq:5.3}
  \tr_X\o(e)^2=0^2\cd d_0 +\dfr{1}{2^2}\, d_{1/2}+\dfr{1}{16^2}\, d_{1/16}
  =\dfr{25-d}{2(5d+22)}.
\end{equation}
The possible integer solution of \eqref{eq:5.2} and \eqref{eq:5.3} are
$(d_0,d_{1/2},d_{1/16})=(3,1,0)$ and $(10,0,8)$.
The case $(d_0,d_{1/2},d_{1/16})=(10,0,8)$ is impossible.
For, if $K$ is isomorphic to $\sqrt{2}A_1$ then $\o(e)$ does not have eigenvalue $1/16$ on $X$,  
so $K$ must be $\sqrt{2}E_8$.
However, in this case we can also find another simple $c=1/2$ Virasoro vector 
$e'\in V_K^+$ such that $\o(e')$ acts on $V_K^-$ semisimply with eigenvalues  
only $0$ and $1/2$. 
(Take a norm 4 element $\alpha \in K$ and consider $V_{\Z \alpha}^+\subset V_K^+$).
But considering $\tr_X \o(e')^2$ we obtain a contradiction.
Therefore, the possible solution is only $(d_0,d_{1/2},d_{1/16})=(3,1,0)$.
We summarize the discussion here in the next theorem.

\begin{thm}\label{thm:5.2}
  Let $L$ be an even positive definite lattice without roots.
  If $L$ contains a norm 4 vector and the weight 1 subspace of $V_L^-$ forms 
  a conformal 4-design then $\mathrm{rank}\, L=\dim (V_L)_1=4$.
\end{thm}

\begin{rem}
  It is shown in \cite{Ma} that if a $c=4$ VOA of OZ-type satisfies: 
  (i)~ it is of class $\mathcal{S}^4$, (ii)~ there is a simple $c=1/2$ Virasoro 
  vector such that its possible eigenvalues of the zero-mode on the Griess algebra 
  are $0$, $1/2$ and $2$, then its Griess algebra is 22-dimensional.
  The fixed point sub VOA $V_{\sqrt{2}D_4}^+$ would be an example with $c=4$ and 
  $\dim V_2=22$ and so $L=\sqrt{2}D_4$ would be an example of a lattice satisfying 
  Theorem \ref{thm:5.2}. 
  (The author expects this is the unique example.)
  $V_{\sqrt{2}D_4}^+$ is isomorphic to the Hamming code VOA and its full automorphism 
  group is $2^6{:}(\mathrm{GL}_3(2)\times \mathrm{S}_3)$ (cf.~\cite{MM,Sh}). 
  The top level of $V_{\sqrt{2}D_4}^-$ is stable under this group.
\end{rem}

\begin{rem}
  Consider the case $L=\sqrt{2}E_8$.
  It is expected (cf.~\cite{Ma}) that $V_{\sqrt{2}E_8}^+$ is of class $\mathcal{S}^6$.
  It is shown in \cite{G,Sh} that $\aut(V_{\sqrt{2}E_8}^+)=\mathrm{O}^+_{10}(2)$.
  The Griess algebra of $V_{\sqrt{2}E_8}^+$ is 156-dimensional which is a direct sum of
  a 1-dimensional module and a 155-dimensional irreducible module over 
  $\mathrm{O}^+_{10}(2)$ (cf.~\cite{ATLAS}).
  It is shown in \cite{G,LSY} that $V_{\sqrt{2}E_8}^+$ has totally 496 simple $c=1/2$ 
  Virasoro vectors which form an $\mathrm{O}^+_{10}(2)$-orbit in the Griess algebra.
  The 240 vectors of them are of $A_1$-type, and the remaining 256 vectors are of 
  $E_8$-type (cf.~\cite{LSY}).
  All $A_1$-types and all $E_8$-types form mutually distinct 
  $2^8{:}\mathrm{O}^+_8(2)$-orbits, where $2^8{:}\mathrm{O}^+_8(2)$ is the 
  centralizer of $\theta$ in $\aut(V_{\sqrt{2}E_8})$.
  The $A_1$-type vector has no eigenvalue 1/16 on $X$ and corresponds to the solution
  $(d_0,d_{1/2},d_{1/16})=(7,1,0)$ of \eqref{eq:5.2}, whereas the $E_8$-type has 
  only eigenvalue 1/16 on $X$ and corresponds to the solution $(d_0,d_{1/2},d_{1/16})=(0,0,8)$
  of \eqref{eq:5.2}.
  It is shown in \cite{Sh} that under the conjugation of modules by automorphisms, 
  the stabilizer of $V_{\sqrt{2}E_8}^-$ in 
  $\aut(V_{\sqrt{2}E_8}^+)=\mathrm{O}^+_{10}(2)$ is $2^8{:}\mathrm{O}^+_8(2)$.
  So $X$ is $2^8{:}\mathrm{O}^+_8(2)$-stable but not $\mathrm{O}^+_{10}(2)$-stable.
\end{rem}

\subsection{The Baby-monster SVOA: the case $h=3/2$}

Here we consider the Baby-monster SVOA $\VB^\natural =\VB^{\natural,0}\oplus \VB^{\natural,1}$ 
introduced by H\"{o}hn in \cite{Ho1} which affords an action of the Baby-monster sporadic finite simple group $\B$.
Let $\theta=(-1)^{2L(0)}\in \aut(\VB^\natural)$ be the canonical $\Z_2$-symmetry.
Then $\VB^\natural$ and $\theta$ satisfy Conditions \ref{cond:1} and \ref{cond:2}.

\begin{rem}
It is shown in \cite{Ho2,Y} that $\aut(\VB^\natural)= \aut(\VB^{\natural,0}) \times \la \theta\ra$, 
$\aut(\VB^{\natural,0})\simeq \B$ and the even part $\VB^{\natural,0}$ has three 
irreducible modules which are all $\B$-stable.
\end{rem}

The top level of $\VB^{\natural,1}$ is of dimension $d=4371$ and has the 
top weight $h=3/2$.
It is shown in Lemma 2.6 of \cite{DLM} (see also \cite{Ho3,Ho4}) that 
$\VB^{\natural,0}$ is of class $\mathcal{S}^6$.
Therefore the top level of $\VB^{\natural.1}$ forms a conformal 6-design based on $\VB^{\natural,0}$.
(Actually, it is shown in \cite{Ho3} and \cite{Ho4} that $\VB^{\natural,0}$ is 
of class $\mathcal{S}^7$ and $\VB^{\natural,1}_{3/2}$ is a conformal 7-design, respectively.)

Let $t$ be a 2A-involution of $\B$.
Then $C_\B(t)\simeq 2{\cd}{}^2E_6(2){:}2$ (cf.~\cite{ATLAS}).
The Griess algebra of $\VB^{\natural,0}$ is of dimension 96256 and 
we have the following decompositions as a module over $\B$ and $C_\B(t)$ (cf.~\cite{HLY}).
\begin{equation}\label{eq:5.4}
\begin{array}{lll}
  \VB^{\natural,0}_2
  &=\ul{\bf{1}}\oplus \ul{\bf{96255}} 
  & \mbox{ over } \B,
  \vsb\\
  &= \ul{\bf{1}}\oplus \ul{\bf{1}}\oplus \ul{\bf{1938}}\oplus \ul{\bf{48620}} \oplus \ul{\bf{45696}}
  & \mbox{ over } C_\B(t).
\end{array}
\end{equation}
The $\B$-invariant subalgebra of the Griess algebra is 1-dimensional 
spanned by the conformal vector $\w$ of $\VB^\natural$, but the $C_\B(t)$-invariant 
subalgebra of the Griess algebra forms a 2-dimensional commutative (and associative) 
subalgebra spanned by two mutually orthogonal Virasoro vectors. 
It is shown in \cite{HLY} that central charges of these Virasoro vectors are 
$7/10$ and $114/5$.
(The sum is the conformal vector of $\VB^\natural$.)
Let $e$ be the shorter one, the simple $c=7/10$ Virasoro vector fixed by $C_\B(t)$.
For the odd part we have the following decompositions.
\begin{equation}\label{eq:5.5}
\begin{array}{lll}
  \VB^{\natural,1}_{3/2}
  &= \ul{\bf 4371} 
  & \mbox{ over } \B,
  \vsb\\
  &= \ul{\bf 1}+\ul{\bf 1938}+\ul{\bf 2432} 
  & \mbox{ over } C_\B(t).
\end{array}
\end{equation}
Since $e$ is fixed by $C_\B(t)$, its zero-mode acts as scalars on 
$C_\B(t)$-irreducible components.
As $c_{4,5}=7/10$ belongs to the minimal series \eqref{eq:3.8}, 
$L(\sfr{7}{10},0)$ has 6 irreducible modules 
$L(\sfr{7}{10},h)$ with $h=0,7/16,3/80,3/2,3/5,1/10$, and $\o(e)$ acts on each 
$C_\B(t)$-irreducible component in \eqref{eq:5.5} by one of these values.
Denote $\lambda_1$, $\lambda_2$ and $\lambda_3$ the eigenvalues of $\o(e)$ on
$\ul{\bf 1}$, $\ul{\bf 1938}$ and $\ul{\bf 2432}$ in \eqref{eq:5.5}, respectively.
Applying Theorem \ref{thm:4.11} we will compute these eigenvalues.
Set $X=\VB^{\natural,1}_{3/2}$.
Since $\VB^{\natural,0}$ is of class $\mathcal{S}^6$ under $\B$ and 
$\VB^{\natural,1}$ is $\B$-stable, by (1), (2), (3) of Theorem \ref{thm:4.11} we have
\begin{equation}\label{eq:5.6}
  \tr_X\o(e)=\dfr{1953}{10},~~
  \tr_X\o(e)^2=\dfr{2163}{100},~~
  \tr_X\o(e)^3=\dfr{5313}{1000}.
\end{equation}
On the other hand, we have
\begin{equation}\label{eq:5.7}
  \tr_X\o(e)^j = \lambda_1^j + 1938 \lambda_2^j +2432 \lambda_3^j
\end{equation}
for $j=1,2,3$.
Solving \eqref{eq:5.6} and \eqref{eq:5.7}, we obtain a unique rational solution
$\lambda_1=3/2$, $\lambda_2=1/10$ and $\lambda_3=0$ which are consistent with 
the representation theory of $L(\sfr{7}{10},0)$.
\begin{equation}\label{eq:5.8}
  \begin{array}{ccccccc}
  \VB^{\natural,1}_{3/2}
  &=& \ul{\bf 1}&+&\ul{\bf 1938}&+&\ul{\bf 2432} 
  \vsb\\
  \o(e)&:& \dfr{3}{2} && \dfr{1}{10} && 0
\end{array}
\end{equation}
Let $x$ be a non-zero vector in the 1-dimensional $C_\B(t)$-invariant subspace of 
$\VB^{\natural,1}_{3/2}$.
Since both $\o(\w)$ and $\o(e)$ act on $x$ by 3/2, it follows $\o(\w-e)x=0$ and 
$x\in \ker (\w-e)_{(0)}$.
This implies $x$ is a square root of $e$ in the extended Griess algebra 
$\VB^{\natural,0}_2\oplus \VB^{\natural,1}_{3/2}$ and $e$ is the conformal vector 
of the subalgebra generated by $x$.
Thus as we discussed in Section \ref{sec:3.1} $\la x\ra$ is isomorphic to 
the $N=1$ $c=7/10$ Virasoro SVOA which is isomorphic to 
$L(\sfr{7}{10},0)\oplus L(\sfr{7}{10},\sfr{3}{2})$ as a $\la e\ra$-module.

\begin{prop}\label{prop:5.6}
  Let $t$ be a 2A-element of $\B$.
  Then $(\VB^\natural)^{C_\B(t)}$ has a full sub SVOA isomorphic to 
  $L(\sfr{114}{5},0)\tensor \l( L(\sfr{7}{10},0)\oplus L(\sfr{7}{10},\sfr{3}{2})\r)$.
\end{prop}

Let us recall the notion of Miyamoto involutions.
A simple $c=7/10$ Virasoro vector $u$ of an SVOA $V$ is called {\it of $\sigma$-type} on $V$ 
if there is no irreducible $\la u\ra\simeq L(\sfr{7}{10},0)$-submodule of $V$ 
isomorphic to either $L(\sfr{7}{10},\sfr{7}{16})$ or $L(\sfr{7}{10},\sfr{3}{80})$.
If $u$ is of $\sigma$-type on $V$, then define a linear automorphism $\sigma_u$ of $V$ 
acting on an irreducible $\la u\ra$-submodule $M$ of $V$ as follows.
\begin{equation}\label{eq:5.9}
  \sigma_u|_M=
  \begin{cases}
    ~\, 1 & \mbox{ if }~ M\simeq L(\sfr{7}{10},0),~ L(\sfr{7}{10},\sfr{3}{5}), 
    \vsb\\
    -1 & \mbox{ if }~ M\simeq L(\sfr{7}{10},\sfr{1}{10}),~ L(\sfr{7}{10},\sfr{3}{2}).
  \end{cases}
\end{equation}
Then $\sigma_u$ is well-defined and the fusion rules of $L(\sfr{7}{10},0)$-modules 
guarantees $\sigma_u\in \aut(V)$ (cf.~\cite{Mi1}). 
It is shown in \cite{HLY} that the map $u \mapsto \sigma_u$ provides 
a one-to-one correspondence between the set of simple $c=7/10$ Virasoro vectors 
of $\VB^{\natural,0}$ of $\sigma$-type and the 2A-conjugacy class of 
$\B=\aut(\VB^{\natural,0})$.
In this correspondence we have to consider only $\sigma$-type $c=7/10$ Virasoro vectors, 
since we also have non $\sigma$-type ones in $\VB^{\natural,0}$.
We can reformulate this correspondence based on the SVOA $\VB^\natural$.
We say a simple $c=7/10$ Virasoro vector $u$ of $\VB^{\natural,0}_2$ is {\it extendable}
if it has a square root $v\in \VB^{\natural,1}_{3/2}$ in the extended Griess algebra
such that $\la v\ra \simeq L(\sfr{7}{10},0)\oplus L(\sfr{7}{10},\sfr{3}{2})$.

Suppose we have an extendable simple $c=7/10$ Virasoro vector 
$u\in \VB^{\natural,0}_2$ and its square root $v\in \VB^{\natural,1}_{3/2}$.
It is shown in \cite{LLY} that the $\Z_2$-graded simple current extension
$L(\sfr{7}{10},0)\oplus L(\sfr{7}{10},\sfr{3}{2})$ has two irreducible 
untwisted modules $L(\sfr{7}{10},0)\oplus L(\sfr{7}{10},\sfr{3}{2})$ and 
$L(\sfr{7}{10},\sfr{1}{10})\oplus L(\sfr{7}{10},\sfr{3}{5})$.
Therefore, $\o(u)$ acts on $X=\VB^{\natural,1}_{3/2}$ semisimply with possible
eigenvalues $0$, $1/10$, $3/5$ and $3/2$.
In particular, $u$ is of $\sigma$-type.
Let $d_\lambda$ be the dimension of $\o(u)$-eigensubspace of $X=\VB^{\natural,1}_{3/2}$ 
with eigenvalue $\lambda$.
Then 
\begin{equation}\label{eq:5.10}
  \tr_X\o(e)^j = 0^j \cd d_0 + \l(\dfr{1}{10}\r)^j \cd d_{1/10}
  +\l(\dfr{3}{5}\r)^j \cd d_{3/5} +\l(\dfr{3}{2}\r)^j \cd d_{3/2}
\end{equation}
for $0\leq j\leq 3$, where we understand $0^0=1$.
By \eqref{eq:5.6} one can solve this linear system and obtain
$d_0=2432$, $d_{1/10}=1938$, $d_{3/5}=0$ and $d_{3/2}=1$, recovering \eqref{eq:5.8}.
(That $d_{3/5}=0$ and $d_{3/2}=1$ can be also shown by the representation theory of 
$L(\sfr{7}{10},0)\oplus L(\sfr{7}{10},\sfr{3}{2})$, see Remark \ref{rem:5.8} below.)
The trace of $\sigma_u$ on $X$ is 
\begin{equation}\label{eq:5.11}
  \tr_X \sigma_u =2432-1-1938=493.
\end{equation}
By \cite{ATLAS} we see that $-\sigma_u$ belongs to the 2A-conjugacy class of $\B$.
Therefore $\sigma_u\theta\in \aut(\VB^\natural)$ is a 2A-element of $\B$ by \cite{Ho2,Y}.
Summarizing, we have the following reformulation of Theorem 5.13 of \cite{HLY}.

\begin{thm}\label{thm:5.7}
There is a one-to-one correspondence between the subalgebras of $\VB^\natural$ 
isomorphic to the $N=1$ $c=7/10$ simple Virasoro SVOA $L(\sfr{7}{10},0)\oplus L(\sfr{7}{10},\sfr{3}{2})$
and the 2A-elements of the Baby-monster $\B$ given by the association 
$u\mapsto \sigma_u \theta$ where $u$ is the conformal vector of the sub SVOA, $\sigma_u$ 
is defined as in \eqref{eq:5.9} and $\theta=(-1)^{2L(0)}$ is the canonical $\Z_2$-symmetry of $\VB^\natural$.
\end{thm}

\begin{rem}\label{rem:5.8}
Let $V=V^0\oplus V^1$ be an SVOA such that $V^1$ has the top weight $3/2$.
Suppose the top level $X=V^1_{3/2}$ forms a conformal 6-design based on $V^0$.
It is shown in \cite{Ho4} that if $\dim X>1$ then the central charge $c$ of $V$ is either 16 or $47/2$.
Suppose further that there is a simple extendable $c=7/10$ Virasoro vector $e$ of $V$.
Then $e$ and its square root generate a sub SVOA $W$ isomorphic to 
$L(\sfr{7}{10},0)\oplus L(\sfr{7}{10},\sfr{3}{2})$.
By the representation theory of $L(\sfr{7}{10},0)\oplus L(\sfr{7}{10},\sfr{3}{2})$ 
(cf.~\cite{LLY}) $V$ is a direct sum of irreducible $W$-submodules isomorphic to 
$L(\sfr{7}{10},0)\oplus L(\sfr{7}{10},\sfr{3}{2})$ or 
$L(\sfr{7}{10},\sfr{1}{10})\oplus L(\sfr{7}{10},\sfr{3}{5})$.
Set $d:=\dim X$ and let $d_\lambda$ be the dimension of $\o(e)$-eigensubspace of $X$ with eigenvalue $\lambda$.
Then $d=d_0+d_{1/10}+d_{3/5}+d_{3/2}$.
We know $d_{3/2}=1$ as $L(\sfr{7}{10},\sfr{3}{2})$ is a simple current 
$L(\sfr{7}{10},0)$-module, and we also have $d_{3/5}=0$ since $X$ is the top level.
Solving \eqref{eq:5.10} in this case by Theorem \ref{thm:4.11} 
we obtain
\begin{equation}\label{eq:5.12}
\begin{array}{l}
  d_{3/5}=-\dfr{7d(2c-47)(10c-7)(82c-37)}{80c(2c-1)(5c+22)(7c+68)},
  \vsb\\
  d_{3/2}=\dfr{d(800c^3-27588c^2+238596c-112133)}{80c(2c-1)(5c+22)(7c+68)}.
\end{array}
\end{equation}
Combining \eqref{eq:5.12} with $d_{3/5}=0$ and $d_{3/2}=1$ we get two possible
solutions $(c,d)=(\sfr{7}{10},1)$ and $(\sfr{47}{2},4371)$.
The case $V=W$ corresponds to the former, and $V=\VB^\natural$ is an
example corresponding to the latter case.
The author expects that $\VB^\natural$ is the unique example of class $\mathcal{S}^6$ 
corresponding to the latter case.
\end{rem}

\section{Appendix}

\subsection{Coefficients in generalized Casimir vectors}\label{app:casimir}


$A^{(2)}_{[2]}=2*h*d$, 
$A^{(3)}_{[3]}=h*d$, 
$A^{(4)}_{[4]}=3*h*d*(c-2*h+4)$, 
$A^{(4)}_{[2,2]}=2*h*(5*h+1)*d$, 
$A^{(5)}_{[5]}=2*h*d*(c-2*h+4)$, 
$A^{(5)}_{[3,2]}=2*h*(5*h+1)*d$, 
$A^{(6)}_{[6]}=4*h*d*(5*c^3+(-15*h+65)*c^2+(-20*h^2-148*h+148)*c-26*h^2+98*h-92)$, 
$A^{(6)}_{[4,2]}=2*h*d*((42*h+8)*c^2+(-84*h^2+349*h+65)*c-134*h^2-86*h-40)$, 
$A^{(6)}_{[3,3]}=(1/2)*h*d*((70*h+15)*c^2+(614*h+136)*c+248*h^2-464*h-64)$, 
$A^{(6)}_{[2,2,2]}=(4/3)*h*d*((70*h^2+42*h+8)*c+29*h^2-57*h-2)$, 
$A^{(7)}_{[7]}=3*h*d*(5*c^3+(-15*h+65)*c^2+(-20*h^2-148*h+148)*c-26*h^2+98*h-92)$, 
$A^{(7)}_{[5,2]}=2*h*d*((28*h+5)*c^2+(-56*h^2+243*h+41)*c-172*h^2-16*h-28)$, 
$A^{(7)}_{[4,3]}=3*h*d*((14*h+3)*c^2+(-28*h^2+106*h+24)*c+38*h^2-70*h-12)$, 
$A^{(7)}_{[3,2,2]}=2*h*d*((70*h^2+42*h+8)*c+29*h^2-57*h-2)$, 
$A^{(8)}_{[8]}=(1/2)*h*d*(350*c^5+(-1260*h+10080)*c^4+(-560*h^2-31735*h+85005)*c^3+(-5040*h^3-17240*h^2-192290*h+194494)*c^2+(-18520*h^3-43840*h^2+20928*h-8184)*c+4344*h^3-32496*h^2+76488*h-57744)$, 
$A^{(8)}_{[6,2]}=2*h*d*((300*h+50)*c^4+(-900*h^2+7312*h+1176)*c^3+(-1200*h^3-18548*h^2+42969*h+6081)*c^2+(-4552*h^3-52960*h^2+32406*h-1466)*c-536*h^3-26880*h^2+5696*h-2808)$, 
$A^{(8)}_{[5,3]}=(1/2)*h*d*((840*h+175)*c^4+(-1680*h^2+19885*h+4188)*c^3+(-25392*h^2+107936*h+23184)*c^2+(-2016*h^3+1832*h^2+4060*h+968)*c+7792*h^3+1776*h^2-31312*h-6816)$, 
$A^{(8)}_{[4,4]}=(3/2)*h*d*((126*h+28)*c^4+(-504*h^2+2787*h+643)*c^3+(504*h^3-7156*h^2+13198*h+3338)*c^2+(3180*h^3+2372*h^2-2480*h+344)*c-2004*h^3+7248*h^2-6036*h-888)$, 
$A^{(8)}_{[4,2,2]}=2*h*d*((630*h^2+366*h+68)*c^3+(-1260*h^3+9159*h^2+4793*h+958)*c^2+(-6942*h^3+11417*h^2-3187*h+210)*c+1114*h^3-654*h^2-3064*h-168)$, 
$A^{(8)}_{[3,3,2]}=h*d*((1050*h^2+645*h+125)*c^3+(16700*h^2+9170*h+1934)*c^2+(3720*h^3+15510*h^2-8662*h+716)*c-1016*h^3+6444*h^2-8692*h-264)$, 
$A^{(8)}_{[2,2,2,2]}=(2/3)*h*d*((1050*h^3+1260*h^2+606*h+108)*c^2+(3305*h^3-498*h^2-701*h+78)*c-251*h^3+918*h^2-829*h-6)$, 
$A^{(9)}_{[9]}=(2/3)*h*d*(210*c^5+(-756*h+6048)*c^4+(-756*h^2-19311*h+50949)*c^3+(-5544*h^3-19676*h^2-120486*h+115622)*c^2+(-23508*h^3-30448*h^2+18372*h-5456)*c+4428*h^3-26232*h^2+52020*h-34536)$, 
$A^{(9)}_{[7,2]}=2*h*d*((225*h+35)*c^4+(-675*h^2+5565*h+826)*c^3+(-900*h^3-14615*h^2+33776*h+4267)*c^2+(-2910*h^3-49778*h^2+29666*h-1378)*c-2350*h^3-23916*h^2+6718*h-2196)$, 
$A^{(9)}_{[6,3]}=4*h*d*((75*h+15)*c^4+(-225*h^2+1747*h+350)*c^3+(-300*h^3-3933*h^2+9193*h+1814)*c^2+(-1642*h^3-3182*h^2+2740*h-88)*c+1814*h^3-2964*h^2-1022*h-612)$, 
$A^{(9)}_{[5,4]}=2*h*d*((126*h+28)*c^4+(-504*h^2+2787*h+643)*c^3+(504*h^3-7156*h^2+13198*h+3338)*c^2+(3180*h^3+2372*h^2-2480*h+344)*c-2004*h^3+7248*h^2-6036*h-888)$, 
$A^{(9)}_{[5,2,2]}=4*h*d*((210*h^2+117*h+21)*c^3+(-420*h^3+3208*h^2+1602*h+302)*c^2+(-3554*h^3+4166*h^2-784*h+40)*c+710*h^3-1032*h^2-746*h-60)$, 
$A^{(9)}_{[4,3,2]}=2*h*d*((630*h^2+381*h+73)*c^3+(-1260*h^3+8694*h^2+4780*h+1010)*c^2+(-3222*h^3+10336*h^2-4022*h+300)*c+98*h^3+1788*h^2-3890*h-156)$, 
$A^{(9)}_{[3,3,3]}=(1/6)*h*d*((1050*h^2+675*h+135)*c^3+(15770*h^2+9144*h+2038)*c^2+(11160*h^3+13348*h^2-10332*h+896)*c-3048*h^3+11328*h^2-10344*h-240)$, 
$A^{(9)}_{[3,2,2,2]}=(4/3)*h*d*((1050*h^3+1260*h^2+606*h+108)*c^2+(3305*h^3-498*h^2-701*h+78)*c-251*h^3+918*h^2-829*h-6)$, 
$A^{(10)}_{[10]}=(6/5)*h*d*(1050*c^6+(-4200*h+52290)*c^5+(-3150*h^2-195019*h+888199)*c^4+(-31500*h^3-160243*h^2-2900235*h+5888368)*c^3+(-33600*h^4-876400*h^3-2224448*h^2-13733560*h+11872408)*c^2+(-189616*h^4-3013900*h^3-3958988*h^2+2767600*h-800016)*c-29792*h^4+816800*h^3-3744448*h^2+6247744*h-3575424)$, 
$A^{(10)}_{[8,2]}=(1/5)*h*d*((19250*h+2800)*c^5+(-69300*h^2+881440*h+124040)*c^4+(-30800*h^3-2898185*h^2+12963179*h+1696856)*c^3+(-277200*h^4-1275240*h^3-35996682*h^2+64729982*h+6705400)*c^2+(-2026600*h^4-3142080*h^3-112130808*h^2+64543216*h-3235248)*c-335864*h^4-8601520*h^3-49036936*h^2+17453488*h-4052928)$, 
$A^{(10)}_{[7,3]}=(3/10)*h*d*((8250*h+1575)*c^5+(-24750*h^2+368615*h+70030)*c^4+(-33000*h^3-978660*h^2+5161264*h+966596)*c^3+(-814640*h^3-10273412*h^2+22761712*h+3992640)*c^2+(-103200*h^4-2524680*h^3-14390328*h^2+9889536*h-487008)*c+139456*h^4+5842880*h^3-11650816*h^2-447872*h-1475328)$, 
$A^{(10)}_{[6,4]}=(12/5)*h*d*((825*h+180)*c^5+(-4125*h^2+35248*h+7832)*c^4+(1650*h^3-149899*h^2+456809*h+104870)*c^3+(6600*h^4+60200*h^3-1328561*h^2+1698137*h+418354)*c^2+(67132*h^4+339830*h^3-202982*h^2-47264*h+18184)*c-33620*h^4-77560*h^3+609380*h^2-579320*h-120480)$, 
$A^{(10)}_{[6,2,2]}=(4/5)*h*d*((8250*h^2+4400*h+760)*c^4+(-24750*h^3+296210*h^2+151096*h+26364)*c^3+(-33000*h^4-801290*h^3+2704347*h^2+1191343*h+213790)*c^2+(-232460*h^4-4589320*h^3+3842818*h^2-286646*h-11132)*c+28644*h^4+857640*h^3-1710804*h^2-212088*h-52032)$, 
$A^{(10)}_{[5,5]}=(1/5)*h*d*((4620*h+1050)*c^5+(-18480*h^2+198719*h+46201)*c^4+(18480*h^3-639607*h^2+2606009*h+632228)*c^3+(443400*h^3-5076574*h^2+9773262*h+2648692)*c^2+(41376*h^4+1074560*h^3+4481200*h^2-3129584*h+288608)*c+293648*h^4-2776160*h^3+7129072*h^2-5166496*h-672384)$, 
$A^{(10)}_{[5,3,2]}=(1/5)*h*d*((46200*h^2+27115*h+5075)*c^4+(-92400*h^3+1636765*h^2+922753*h+178072)*c^3+(-2371200*h^3+14798146*h^2+7158654*h+1515620)*c^2+(-110880*h^4-9063800*h^3+17179184*h^2-5440168*h+338704)*c-319088*h^4+2238560*h^3-652432*h^2-4969184*h-254976)$, 
$A^{(10)}_{[4,4,2]}=(3/5)*h*d*((6930*h^2+4180*h+800)*c^4+(-27720*h^3+226185*h^2+133283*h+26662)*c^3+(27720*h^4-631060*h^3+1836386*h^2+928514*h+211140)*c^2+(337140*h^4-1677780*h^3+2472564*h^2-812508*h+53544)*c+100532*h^4-226520*h^3+432268*h^2-745384*h-34656)$, 
$A^{(10)}_{[4,3,3]}=(3/10)*h*d*((11550*h^2+7425*h+1485)*c^4+(-23100*h^3+389020*h^2+243566*h+51294)*c^3+(-426920*h^3+3236462*h^2+1740898*h+429980)*c^2+(-81840*h^4+1025620*h^3+3455108*h^2-2154976*h+167248)*c-17296*h^4-522560*h^3+2125936*h^2-1996768*h-52992)$, 
$A^{(10)}_{[4,2,2,2]}=4*h*d*((2310*h^3+2706*h^2+1276*h+224)*c^3+(-4620*h^4+48797*h^3+50252*h^2+22925*h+4434)*c^2+(-47038*h^4+140169*h^3-6264*h^2-27525*h+2578)*c-4966*h^4+9340*h^3+15382*h^2-28252*h-288)$, 
$A^{(10)}_{[3,3,2,2]}=h*d*((11550*h^3+14025*h^2+6809*h+1222)*c^3+(274840*h^3+284503*h^2+133429*h+26346)*c^2+(40920*h^4+764986*h^3-116882*h^2-177916*h+18992)*c+8648*h^4-92288*h^3+251080*h^2-201520*h-1344)$, 
$A^{(10)}_{[2,2,2,2,2]}=(4/15)*h*d*((11550*h^4+23100*h^3+20130*h^2+8580*h+1440)*c^2+(76675*h^4+30590*h^3-25615*h^2-10898*h+1608)*c+3767*h^4-18410*h^3+29929*h^2-16342*h-24)$.

\subsection{Coefficients in the trace formulae}\label{app:formula}

$\Sym (a^0|\w)(a^1|a^2)=(a^0|\w)(a^1|a^2)+(a^1|\w)(a^0|a^2)+(a^2|\w)(a^0|a^1)$,
\\
$\Sym (a^0|\w)(a^1|\w)(a^2|a^3)= (a^0|\w)(a^1|\w)(a^2|a^3)+(a^0|\w)(a^2|\w)(a^1|a^3)
+(a^0|\w)(a^3|\w)(a^1|a^2)+(a^1|\w)(a^2|\w)(a^0|a^3)+(a^1|\w)(a^3|\w)(a^0|a^2)+(a^2|\w)(a^3|\w)(a^0|a^1)$, 
\\
$\Sym (a^0|\w)(a^1|a^2|a^3)=(a^0|\w)(a^1|a^2|a^3)+(a^1|\w)(a^0|a^2|a^3)
+(a^2|\w)(a^0|a^1|a^3)+(a^3|\w)(a^0|a^1|a^2)$, 
\\
$\Sym (a^0|a^1)(a^2|a^3)=(a^0|a^1)(a^2|a^3)+(a^0|a^2)(a^1|a^3)+(a^0|a^3)(a^1|a^2)$, 
\\
$\Sym (a^0|\w)(a^1|\w)(a^2|\w)(a^3|a^4)
=(a^0|\w)(a^1|\w)(a^2|\w)(a^3|a^4)
+(a^0|\w)(a^1|\w)(a^3|\w)(a^2|a^4)
+(a^0|\w)(a^1|\w)(a^4|\w)(a^2|a^3)
+(a^0|\w)(a^2|\w)(a^3|\w)(a^1|a^4)
+(a^0|\w)(a^2|\w)(a^4|\w)(a^1|a^3)
+(a^0|\w)(a^3|\w)(a^4|\w)(a^1|a^2)
+(a^1|\w)(a^2|\w)(a^3|\w)(a^0|a^4)
+(a^1|\w)(a^2|\w)(a^4|\w)(a^0|a^3)
+(a^1|\w)(a^3|\w)(a^4|\w)(a^0|a^2)
+(a^2|\w)(a^3|\w)(a^4|\w)(a^0|a^1)$, 
\\
$\Sym (a^0|\w)(a^1|\w)(a^2|a^3|a^4)
=(a^0|\w)(a^1|\w)(a^2|a^3|a^4)
+(a^0|\w)(a^2|\w)(a^1|a^3|a^4)
+(a^0|\w)(a^3|\w)(a^1|a^2|a^4)
+(a^0|\w)(a^4|\w)(a^1|a^2|a^3)
+(a^1|\w)(a^2|\w)(a^0|a^3|a^4)
+(a^1|\w)(a^3|\w)(a^0|a^2|a^4)
+(a^1|\w)(a^4|\w)(a^0|a^2|a^3)
+(a^2|\w)(a^3|\w)(a^0|a^1|a^4)
+(a^2|\w)(a^4|\w)(a^0|a^1|a^3)
+(a^3|\w)(a^4|\w)(a^0|a^1|a^2)$, 
\\
$\Sym (a^0|\w)(a^1|a^2)(a^3|a^4)
=(a^0|\w)(a^1|a^2)(a^3|a^4)
+(a^0|\w)(a^1|a^3)(a^2|a^4)
+(a^0|\w)(a^1|a^4)(a^2|a^3)
+(a^1|\w)(a^0|a^2)(a^3|a^4)
+(a^1|\w)(a^0|a^3)(a^2|a^4)
+(a^1|\w)(a^0|a^4)(a^2|a^3)
+(a^2|\w)(a^0|a^1)(a^3|a^4)
+(a^2|\w)(a^0|a^3)(a^1|a^4)
+(a^2|\w)(a^0|a^4)(a^1|a^3)
+(a^3|\w)(a^0|a^1)(a^2|a^4)
+(a^3|\w)(a^0|a^2)(a^1|a^4)
+(a^3|\w)(a^0|a^4)(a^1|a^2)
+(a^4|\w)(a^0|a^1)(a^2|a^3)
+(a^4|\w)(a^0|a^2)(a^1|a^3)
+(a^4|\w)(a^0|a^3)(a^1|a^2)$, 
\\
$\Sym (a^0|a^1)(a^2|a^3|a^4)
=(a^0|a^1)(a^2|a^3|a^4)
+(a^0|a^2)(a^1|a^3|a^4)
+(a^0|a^3)(a^1|a^2|a^4)
+(a^0|a^4)(a^1|a^2|a^3)
+(a^1|a^2)(a^0|a^3|a^4)
+(a^1|a^3)(a^0|a^2|a^4)
+(a^1|a^4)(a^0|a^2|a^3)
+(a^2|a^3)(a^0|a^1|a^4)
+(a^2|a^4)(a^0|a^1|a^3)
+(a^3|a^4)(a^0|a^1|a^2)$.

$F^{(3)}_0 = 8*h*d*((70*h^2+42*h+8)*c+29*h^2-57*h-2)$, 
$F^{(3)}_1 = -4*h*d*((14*h+4)*c^2+(-308*h^2-93*h-1)*c+170*h^2+34*h)$, 
$F^{(3)}_2 = h*d*(4*c^3+(-222*h-1)*c^2+(3008*h^2+102*h)*c-1496*h^2)$, 
$F^{(4)}_0 = 16*h*d*((1050*h^3+1260*h^2+606*h+108)*c^2+(3305*h^3-498*h^2-701*h+78)*c-251*h^3+918*h^2-829*h-6)$, 
$F^{(4)}_1 = -8*h*d*((210*h^2+162*h+36)*c^3+(-4620*h^3-3227*h^2-861*h+26)*c^2+(-5614*h^3+2915*h^2-485*h-2)*c-1334*h^3+2622*h^2+92*h)$, 
$F^{(4)}_2 = 2*h*d*(60*h*c^4+(-3330*h^2-523*h-487)*c^3+(45120*h^3+9648*h^2+13856*h-2336)*c^2+(36376*h^3-91186*h^2+43550*h-2232)*c-6760*h^3-47796*h^2+19756*h-696)$, 
$F^{(4)}_3 = 4*h*d*((42*h+36)*c^4+(-1848*h^2-1279*h+513)*c^3+(20328*h^3+13052*h^2-14654*h+2334)*c^2+(-35836*h^3+98516*h^2-43320*h+2232)*c-16700*h^3+43104*h^2-19756*h+696)$, 
$F^{(4)}_4 = (1/2)*h*d*((1128*h+199)*c^4+(-46392*h^2-3311*h-1768)*c^3+(497472*h^3-19488*h^2+73544*h-16440)*c^2+(351008*h^3-726256*h^2+326804*h-17160)*c-72848*h^3-344832*h^2+158048*h-5568)$, 
$F^{(4)}_5 = (-1/2)*h*d*(60*c^5+(-2976*h+1023)*c^4+(44184*h^2-41669*h+2850)*c^3+(-164544*h^3+426432*h^2-65116*h-716)*c^2+(-22112*h^3+23984*h^2+13092*h-1528)*c+68816*h^3-150144*h^2+25024*h)$, 
$F^{(4)}_6 = (1/2)*h*d*(60*c^5+(-2640*h+1311)*c^4+(29400*h^2-51901*h+6954)*c^3+(-1920*h^3+530848*h^2-182348*h+17956)*c^2+(-308800*h^3+812112*h^2-333468*h+16328)*c-64784*h^3+194688*h^2-133024*h+5568)$, 
$F^{(5)}_0=  32*h*d*((11550*h^4+23100*h^3+20130*h^2+8580*h+1440)*c^2+(76675*h^4+30590*h^3-25615*h^2-10898*h+1608)*c+3767*h^4-18410*h^3+29929*h^2-16342*h-24)$, 
$F^{(5)}_1=  -16*h*d*((2310*h^3+3366*h^2+1848*h+360)*c^3+(-50820*h^4-64063*h^3-39624*h^2-9203*h+402)*c^2+(-190058*h^4+21757*h^3+50420*h^2-8593*h-6)*c+14558*h^4-53244*h^3+48082*h^2+348*h)$, 
$F^{(5)}_2=  (4/5)*h*d*((3300*h^2+660*h-40)*c^4+(-183150*h^3-90835*h^2-94567*h-25578)*c^3+(2481600*h^4+1334700*h^3+2540131*h^2+285789*h-163830)*c^2+(7115560*h^4-13778670*h^3+2299334*h^2+2630452*h-245456)*c+858872*h^4+1045920*h^3-6623912*h^2+2211696*h-37056)$, 
$F^{(5)}_3=  (8/5)*h*d*((2310*h^2+3300*h+920)*c^4+(-101640*h^3-123925*h^2+2681*h+13794)*c^3+(1118040*h^4+1178540*h^3-631298*h^2-179402*h+81900)*c^2+(-228580*h^4+4993420*h^3-750692*h^2-1313196*h+122728)*c+344284*h^4-2043720*h^3+3258596*h^2-1105848*h+18528)$, 
$F^{(5)}_4=  (1/5)*h*d*(500*c^5+(62040*h^2+7735*h+25115)*c^4+(-2551560*h^3-564175*h^2-1452063*h+94428)*c^3+(27360960*h^4+2744160*h^3+30534534*h^2-6099454*h-348380)*c^2+(64210400*h^4-208744320*h^3+91532216*h^2+3799848*h-935504)*c-909872*h^4-64093920*h^3+8306672*h^2+7635744*h-148224)$, 
$F^{(5)}_5=  (-1/5)*h*d*((3300*h+1500)*c^5+(-163680*h^2-585*h+39235)*c^4+(2430120*h^3-2017145*h^2-1433609*h+240084)*c^3+(-9049920*h^4+31066560*h^3+13487402*h^2-7601082*h+456180)*c^2+(-41190560*h^4-16962080*h^3+49902728*h^2-9779816*h-107152)*c-808336*h^4+9987680*h^3-17678384*h^2+1778272*h+188928)$, 
$F^{(5)}_6=  (1/5)*h*d*((3300*h+1500)*c^5+(-145200*h^2+25815*h+46595)*c^4+(1617000*h^3-3008545*h^2-1412161*h+350436)*c^3+(-105600*h^4+40494880*h^3+8437018*h^2-9036298*h+1111380)*c^2+(-43019200*h^4+22985280*h^3+43897192*h^2-20285384*h+874672)*c+1945936*h^4-6362080*h^3+8390384*h^2-7068512*h+337152)$, 
$F^{(5)}_7=  (-2/5)*h*d*((660*h+460)*c^5+(-51150*h^2-33647*h+6897)*c^4+(1302180*h^3+829156*h^2-451426*h+40950)*c^3+(-10919040*h^4-7782640*h^3+9315274*h^2-2071698*h+61364)*c^2+(21416272*h^4-59346500*h^3+27298188*h^2-1866624*h+9264)*c+9686000*h^4-25000320*h^3+11458480*h^2-403680*h)$, 
$F^{(5)}_8=  (-1/2)*h*d*(100*c^6+(-8078*h+2861)*c^5+(221174*h^2-203081*h+19684)*c^4+(-2214880*h^3+4802538*h^2-965274*h+52252)*c^3+(4236288*h^4-38346896*h^3+13282628*h^2-1695920*h+25584)*c^2+(12825792*h^4-32289856*h^3+12276272*h^2-758816*h+17536)*c-155904*h^4-1722368*h^3+4176000*h^2-215296*h)$, 
$F^{(5)}_{01423}= (-1/10)*h*d*(500*c^6+(-33130*h+25625)*c^5+(707230*h^2-1434751*h+485426)*c^4+(-4128560*h^3+27827338*h^2-20020070*h+4414912)*c^3+(-15989760*h^4-192127280*h^3+258582588*h^2-134080200*h+15339472)*c^2+(136946816*h^4-932999600*h^3+1000619648*h^2-291521120*h+11640256)*c+25836352*h^4-152808960*h^3+252283328*h^2-102259584*h+4475904)$, 
$F^{(5)}_{01324}= (1/10)*h*d*(500*c^6+(-35770*h+23785)*c^5+(911830*h^2-1300163*h+457838)*c^4+(-9337280*h^3+24510714*h^2-18214366*h+4251112)*c^3+(27686400*h^4-160996720*h^3+221321492*h^2-125793408*h+15094016)*c^2+(51281728*h^4-695613600*h^3+891426896*h^2-284054624*h+11603200)*c-12907648*h^4-52807680*h^3+206449408*h^2-100644864*h+4475904)$, 
$F^{(5)}_{12034}= (1/10)*h*d*(100*c^6+(-3150*h+5575)*c^5+(15650*h^2-148721*h+119806)*c^4+(-550800*h^3-1490922*h^2-4041146*h+949728)*c^3+(14745600*h^4+53833840*h^3+29411876*h^2-25511768*h+3284592)*c^2+(-167754624*h^4+73524400*h^3+137639360*h^2-61086944*h+2624448)*c+3115968*h^4-16839680*h^3+31808832*h^2-21690496*h+1007616)$, 
$F^{(5)}_{01234}=  h*d*(50*c^6+(-5304*h+1238)*c^5+(204604*h^2-99615*h+13827)*c^4+(-3383208*h^3+2721160*h^2-892294*h+59250)*c^3+(20120832*h^4-26868960*h^3+17636364*h^2-3522876*h-65568)*c^2+(41237472*h^4-107859720*h^3+46229896*h^2-423632*h-256192)*c-772320*h^4-27659904*h^3+9574560*h^2+1892928*h-31488)$, 
$F^{(5)}_{02134}= (-1/10)*h*d*(100*c^6+(150*h+7215)*c^5+(-237790*h^2-280007*h+130072)*c^4+(5835360*h^3+1550206*h^2-5491014*h+812348)*c^3+(-38223360*h^4+30779440*h^3+63018188*h^2-27827312*h+1723824)*c^2+(-139274688*h^4-126848480*h^3+208934224*h^2-38241856*h-180160)*c+69888*h^4+19482880*h^3-41267328*h^2+3697664*h+692736)$, 
$F^{(5)}_{03214}=F^{(5)}_{04213}= (1/5)*h*d*((1650*h+820)*c^5+(-126720*h^2-65643*h+5133)*c^4+(3193080*h^3+1520564*h^2-724934*h-68690)*c^3+(-26484480*h^4-11527200*h^3+16803156*h^2-1157772*h-780384)*c^2+(14239968*h^4-100186440*h^3+35647432*h^2+11422544*h-1402304)*c-1523040*h^4+18161280*h^3-36538080*h^2+12694080*h-157440)$, 
$F^{(5)}_{02413}=F^{(5)}_{03412}= (-1/5)*h*d*((3630*h+5660)*c^5+(-199320*h^2-209673*h+193503)*c^4+(3472920*h^3+1907324*h^2-7596850*h+2076826)*c^3+(-18585600*h^4-196400*h^3+96084724*h^2-62800300*h+7605776)*c^2+(36408928*h^4-385775160*h^3+469619144*h^2-143863520*h+5776288)*c+13307936*h^4-72098560*h^3+115701664*h^2-50591552*h+2237952)$, 
$F^{(5)}_{02314}=F^{(5)}_{04312}= (1/5)*h*d*((2310*h+4740)*c^5+(-97020*h^2-142379*h+179709)*c^4+(868560*h^3+249012*h^2-6693998*h+1994926)*c^3+(3252480*h^4+15368880*h^3+77454176*h^2-58656904*h+7483048)*c^2+(-6423616*h^4-267082160*h^3+415022768*h^2-140130272*h+5757760)*c-6064064*h^4-22097920*h^3+92784704*h^2-49784192*h+2237952)$, 
$F^{(5)}_{04123}=F^{(5)}_{03124}=(-2/5)*h*d*(150*c^6+(-10060*h+5380)*c^5+(217020*h^2-323853*h+57123)*c^4+(-1309760*h^3+6390724*h^2-2579346*h+268402)*c^3+(-4260480*h^4-40238760*h^3+32357832*h^2-9076728*h+462936)*c^2+(-18786432*h^4-72074440*h^3+67578896*h^2-10508984*h-23120)*c-177408*h^4+2717760*h^3-5096832*h^2+655296*h+173184)$, 
$F^{(5)}_{14023}=F^{(5)}_{13024}= (1/5)*h*d*(300*c^6+(-21770*h+9940)*c^5+(560760*h^2-582063*h+109113)*c^4+(-5812600*h^3+11260884*h^2-4433758*h+605494)*c^3+(17963520*h^4-68950320*h^3+47912508*h^2-16995684*h+1706256)*c^2+(-51812832*h^4-43962440*h^3+99510360*h^2-32440512*h+1356064)*c+1168224*h^4-12725760*h^3+26344416*h^2-11383488*h+503808)$, 
$E^{(3)}_3= 8*h*d*((70*h^2+42*h+8)*c+29*h^2-57*h-2)$, 
$E^{(3)}_2= -12*h*d*((14*h+4)*c^2+(-308*h^2-93*h-1)*c+170*h^2+34*h)$, 
$E^{(3)}_1= 2*h*d*(4*c^3+(-222*h-1)*c^2+(3008*h^2+102*h)*c-1496*h^2)$, 
$E^{(4)}_4= 16*h*d*((1050*h^3+1260*h^2+606*h+108)*c^2+(3305*h^3-498*h^2-701*h+78)*c-251*h^3+918*h^2-829*h-6)$, 
$E^{(4)}_3= -48*h*d*((210*h^2+162*h+36)*c^3+(-4620*h^3-3227*h^2-861*h+26)*c^2+(-5614*h^3+2915*h^2-485*h-2)*c-1334*h^3+2622*h^2+92*h)$, 
$E^{(4)}_2= 4*h*d*((366*h+108)*c^4+(-18864*h^2-5929*h-409)*c^3+(241464*h^3+77748*h^2+11462*h-2342)*c^2+(37996*h^3-69196*h^2+44240*h-2232)*c-77140*h^3-61872*h^2+19756*h-696)$, 
$E^{(4)}_1= 2*h*d*((1464*h+487)*c^4+(-61176*h^2-13543*h+2336)*c^3+(660096*h^3+84928*h^2-43688*h+2232)*c^2+(64320*h^3+61872*h^2-19756*h+696)*c-206448*h^3)$, 
$E^{(5)}_5= 32*h*d*((11550*h^4+23100*h^3+20130*h^2+8580*h+1440)*c^2+(76675*h^4+30590*h^3-25615*h^2-10898*h+1608)*c+3767*h^4-18410*h^3+29929*h^2-16342*h-24)$, 
$E^{(5)}_4= -160*h*d*((2310*h^3+3366*h^2+1848*h+360)*c^3+(-50820*h^4-64063*h^3-39624*h^2-9203*h+402)*c^2+(-190058*h^4+21757*h^3+50420*h^2-8593*h-6)*c+14558*h^4-53244*h^3+48082*h^2+348*h)$, 
$E^{(5)}_3= 8*h*d*((13530*h^2+11220*h+2680)*c^4+(-671220*h^3-553445*h^2-181091*h-9774)*c^3+(8317320*h^4+6205020*h^3+3186368*h^2+33372*h-81960)*c^2+(13545380*h^4-12577080*h^3+2346592*h^2+1321316*h-122728)*c+2750596*h^4-4039320*h^3-3472036*h^2+1105848*h-18528)$, 
$E^{(5)}_2= -4*h*d*((1320*h+420)*c^5+(-182820*h^2-101429*h-18681)*c^4+(5969040*h^3+3213887*h^2+527763*h-122880)*c^3+(-58143360*h^4-27737760*h^3-6853602*h^2+3391274*h-184092)*c^2+(-19549216*h^4+50103960*h^3-30930304*h^2+2972472*h-27792)*c+17527600*h^4+30443040*h^3-11458480*h^2+403680*h)$, 
$E^{(5)}_1= -2*h*d*(100*c^6+(1470*h+6495)*c^5+(-501790*h^2-424803*h+40956)*c^4+(15693120*h^3+8719374*h^2-2073438*h+61364)*c^3+(-141373440*h^4-54143280*h^3+27458268*h^2-1866624*h+9264)*c^2+(-12282432*h^4-30443040*h^3+11458480*h^2-403680*h)*c+47895936*h^4)$.


\end{document}